\RequirePackage{tikz}
\documentclass[pdflatex,sn-mathphys]{sn-jnl}

\usepackage{enumitem}
\usepackage{pgfplots}
\usepackage{tikz}
\usepackage{subcaption}
\usepackage{mathtools}
\pgfplotsset{compat=1.18}

\jyear{2021}%

\theoremstyle{thmstyleone}%
\theoremstyle{thmstyletwo}%

\theoremstyle{thmstylethree}%

\newtagform{bluetag}{\color{blue}(}{)}

\begin{document}

\title[Learning-Assisted Optimization for Transmission Switching]{Learning-Assisted Optimization for Transmission Switching}


\author*[1,3]{\fnm{Salvador} \sur{Pineda}}\email{spineda@uma.es}

\author[2,3]{\fnm{Juan Miguel} \sur{Morales}}\email{juan.morales@uma.es}

\author[2,3]{\fnm{Asunción} \sur{Jiménez-Cordero}}\email{asuncionjc@uma.es}

\affil*[1]{\orgdiv{Dept. of Electrical Engineering}, \orgname{University of Málaga},  \country{Spain}}

\affil[2]{\orgdiv{Dept. of Mathematical Analysis, Statistics and Operations Research, and Applied Mathematics}, \orgname{University of Málaga}, \country{Spain}}

\affil[3]{\orgdiv{OASYS Research Group}, \orgname{University of Málaga},  \country{Spain}}


\abstract{The design of new strategies that exploit methods from Machine Learning to facilitate the resolution of challenging and large-scale mathematical optimization problems has recently become an avenue of prolific and promising research. In this paper, we propose a novel learning procedure to assist in the solution of a well-known computationally difficult optimization problem in power systems: The Direct Current Optimal Transmission Switching (DC-OTS)  problem. The DC-OTS problem consists in finding the configuration of the power network that results in the cheapest dispatch of the power generating units.  With the increasing variability in the operating conditions of power grids, the DC-OTS problem has lately sparked renewed interest, because operational strategies that include topological network changes have proved to be effective and efficient in helping maintain the balance between generation and demand.  The DC-OTS problem includes a set of binaries that determine the on/off status of the switchable transmission lines. Therefore, it takes the form of a mixed-integer program, which is NP-hard in general. In this paper, we propose an approach to tackle the DC-OTS problem that leverages known solutions to past instances of the problem to speed up the mixed-integer optimization of a new unseen model. Although our approach does not offer optimality guarantees,} a series of numerical experiments run on a real-life power system dataset show that it features a very high success rate in identifying the optimal grid topology (especially when compared to alternative competing heuristics), while rendering remarkable speed-up factors.




\maketitle

\section{Introduction}\label{sec:intro}

Power systems are colossal and complex networks engineered to reliably supply electricity where it is needed at the \emph{lowest} possible cost. For this, operational routines based on the \emph{Optimal Power Flow} (OPF) problem are executed daily and in real time to guarantee the most cost-efficient dispatch of power generating units that satisfy the grid constraints. In particular, the way power flows through a power network is determined by the so-called \emph{Kirchhoff's laws}. These laws are responsible for the fact that \emph{switching off} a transmission line in the grid can actually result in a lower electricity production cost (a type of ``Braess' Paradox'') and have provided power system operators with a complementary control action, namely, changes in the grid \emph{topology}, to reduce this cost even further. The possibility of flexibly exploiting the topological configuration of the grid was first suggested in~\cite{o2005dispatchable} and later formalized in~\cite{fisher2008optimal} into what we know today as the \emph{Optimal Transmission Switching} (OTS) problem. Essentially, the OTS problem is the OPF problem enriched with a whole new set of on/off variables that model the status of each \emph{switchable} transmission line in the system. The OPF formulation we use as a basis to pose the OTS problem is built on the widely used \emph{direct current} (DC) linear approximation of the power flow equations. Even so, the resulting formulation of the OTS problem, known as DC-OTS, takes the form of a mixed-integer program, which has been proven to be NP-hard for general network classes~\cite{kocuk2016cycle, fattahi2019bound}.

Thus, the DC-OTS problem consists in finding the configuration of the power network that results in the cheapest dispatch of the power generating units subject to constraints such as thermal limits on transmission lines, generating units' capacity bounds, and network connectivity conditions. To date, the resolution of the DC-OTS has been approached from two distinct methodological points of view, namely, by means of \emph{exact} methods and by way of \emph{heuristics}. The former exploit techniques from mixed-integer programming such as bounding, tightening, and the generation of valid cuts to solve the DC-OTS to (certified) global optimality, while the latter seek to quickly identify good solutions of the problem, but potentially forgoing optimality and even at the risk of suggesting infeasible grid topologies. 

Among the methods that are exact, we highlight the works in~\cite{kocuk2016cycle}, \cite{fattahi2019bound}, \cite{ruiz2016security}, and~\cite{dey2022node}. More specifically, the authors in~\cite{kocuk2016cycle} propose a cycle-based formulation of the DC-OTS problem, which results in a mixed-integer linear program. They prove the NP-hardness of the DC-OTS even if the power grid takes the form of a series-parallel graph with only one generation-demand pair, and derive classes of strong valid inequalities for a relaxation of their formulation that can be separated in polynomial time. In \cite{fattahi2019bound}, the authors work instead with the mixed-integer linear formulation of the DC-OTS that employs a big-M to model the disjunctive nature of the equation linking the power flow through a switchable line and the voltage angles at the nodes the line connects. This is the formulation of the DC-OTS we also consider in this paper. The big-M must be a valid upper bound of the maximum angle difference when the switchable line is open. In \cite{fattahi2019bound}, it is proven that determining this maximum is NP-hard and, consequently, they propose to set the big-M to the shortest path between the nodes concerned over a spanning subgraph that is assumed to exist. The authors in~\cite{ruiz2016security} conduct a computational study of a mixed-integer linear reformulation of the DC-OTS problem alternative to that considered in \cite{fattahi2019bound}. This reformulation makes use of the so-called \emph{power transfer distribution factors} (PTDFs) and the notion of \emph{flow-cancelling transactions} to model open lines. They argue that this reformulation comparatively offers significant computational advantages, especially for large systems and when the number of switchable lines is relatively small.
Finally, a family of cutting planes for the DC-OTS problem are developed in~\cite{dey2022node}. These cutting planes are derived from the polyhedral description of the integer hull of a certain constraint set that appears in the DC-OTS problem. Specifically, this constraint set is made up of a nodal power balance equation together with the power flow limits of the associated incident lines. Those of these limits that correspond to switchable lines are multiplied by the respective binary variable. 

In practice, though, the complexity and size of real-life power grids often render exact solutions computationally infeasible. Therefore, heuristics, or approximate solution methods, become essential for tackling the DC-OTS efficiently. Among the heuristic methods that have been proposed in the technical literature, we can distinguish two main groups. The first group includes the heuristic approaches that do not rely on the solutions of previous instances of the OTS problem. For example, some heuristics trim down the computational time by reducing the number of lines that can be switched off \cite{liu2012heuristic, barrows2012computationally, flores2020alternative}. While these approaches do not reach the maximum cost savings, the reported numerical studies show that the cost increase with respect to the optimal solution is small in most cases. Other related works maintain the original set of switchable lines and determine their on/off status using greedy algorithms \cite{fuller2012fast, crozier2022feasible}. They use dual information of the OPF problem to rank the lines according to the impact of their status on the operational cost. Finally, the authors of \cite{hinneck2022optimal} propose solving the OTS problem in parallel with heuristics that generate good candidate solutions to speed up conventional MIP algorithms. The second group comprises  data-based heuristic methods that require information about the optimal solution of past OTS problems. For instance, the authors of \cite{johnson2020knearest} use a $K$-nearest neighbor strategy to drastically reduce the search space of the integer solution to the DC-OTS problem. In particular, given a collection of past instances of the problem (whose solution is assumed to be known and available), they restrict the search space to the $K$ integer solutions of those instances which are the closest to the one to be solved in terms of the problem parameters (for example, nodal demands). They then provide as solution to the instance of the DC-OTS problem under consideration the one that results in the lowest cost. This last step requires solving $K$ linear programs, one per candidate integer solution. Conversely, various alternative data-driven methods, distinct from the $K$ nearest neighbor, have also been explored to enhance the solution of the DC-OTS problem. For example, references \cite{yang2019line, han2022learning, bugaje2023real} present sophisticated methodologies to learn the status of switchable lines using neural networks. 

Against this background, in this paper, we propose a novel method to address the DC-OTS  by exploiting known solutions to \emph{past} instances of the problem. Indeed, according to \cite{bengio2021machine, parmentier2022learning}, our approach aligns with machine learning strategies that extract valuable insights from prior solutions of an optimization problem, subsequently applying this knowledge to address new, unseen instances. Specifically, our approach leverages information from previous instances in two distinct yet potentially synergistic ways. First, from these past solutions, we infer those switchable lines that are most likely to be operational (resp. inoperative) in the current instance of the problem (the one we want to solve). Mathematically, this translates into fixing a few binaries to one (resp. zero), an apparently small action that brings, however, substantial benefits in terms of computational speed. Second, beyond the speed-up that one can expect from simply reducing the number of binaries in a MILP, this strategy also allows us to leverage the shortest-path-based argument invoked in \cite{fattahi2019bound} to further tighten the big-Ms in the problem formulation, with the consequent extra computational gain. 

Alternatively, we also investigate the potential of directly \emph{inferring} the big-M values from past solutions to the problem, eliminating the need for the shortest-path calculation. In any case, the inference of the binaries to be fixed and/or the values of the big-Ms to be used is conducted through a Machine Learning algorithm of the decision-maker's choice. In this paper, we have opted for the use of the $K$-nearest neighbors methodology due to its simplicity, as well as its interpretability and low computational time required for the training task. Besides, this approach has demonstrated success in mitigating the complexity of related challenges, such as the widely studied DC Unit Commitment problem, as evidenced by prior works \cite{pineda2020data, jimenez2022warm}.

Importantly, while our proposal is not endowed with theoretical guarantees of optimality (and thus, belongs to the group of heuristics discussed above), the role that Machine Learning plays in it is supportive rather than surrogative (we still need to solve the MILP problem), which results in significantly lower rates of infeasibility and suboptimality, as demonstrated in the numerical experiments.

The remainder of this paper is structured as follows. Section~\ref{sec:ots} introduces the DC-OTS problem mathematically and discusses how to equivalently reformulate it as a mixed-integer linear program (MILP) through the use of large enough constants (the so-called \emph{big-Ms}). Section~\ref{sec:met} describes the different methods we consider in this paper to identify the most cost-efficient grid topology of a power system, including those we propose and those we use for benchmarking. A series of numerical experiments run on a 118-bus power system typically used in the context of the DC-OTS problem are presented and discussed in Section~\ref{sec: num}. Finally, conclusions and further research are duly drawn in Section~\ref{sec:concl}.

\section{Optimal transmission switching}\label{sec:ots}
We start this section by introducing the standard and well-known formulation of the \emph{Direct Current Optimal Transmission Switching} problem (DC-OTS), which will serve us a basis to construct and motivate its mixed-integer reformulation immediately after. 

Consider a power network consisting of a collection of nodes $\mathcal{N}$ and transmission lines $\mathcal{L}$. To lighten the mathematical formulation of the DC-OTS, we assume \emph{w.l.o.g} that there is one generator and one power load per node $n \in \mathcal{N}$. The power dispatch of the generator and the power consumed by the power load are denoted by $p_n$ and $d_n$, respectively. Each generator is characterized by a minimum and maximum power output, $\underline{p}_n$ and $\overline{p}_n$, and a marginal production cost $c_n$. We represent the power flow through the line $(n,m) \in \mathcal{L}$ connecting nodes $n$ and $m$ by $f_{nm}$, with $f_{nm} \in [-\overline{f}_{nm}, \overline{f}_{nm}]$. For each node $n$ we distinguish between the set of transmission lines whose power flow \emph{enters} the node, $\mathcal{L}_{n}^{+}$, and the set of transmission lines whose power flow \emph{leaves} it, $\mathcal{L}_{n}^{-}$. The power network includes a subset $\mathcal{L}_\mathcal{S} \subseteq \mathcal{L}$ of lines that can be switched on/off. If the line $(n,m) \in \mathcal{L}_\mathcal{S}$, its status is determined by a binary variable $x_{nm}$, which takes value 1 if the line is fully operational, and 0 when disconnected. In a DC power network, the flow $f_{nm}$ through an operational line is given by the product of the susceptance of the line, $b_{nm}$, and the difference of the voltage angles at nodes $n$ and $m$, i.e., $\theta_n-\theta_m$. We use bold symbols to define the vectors of variables $\mathbf{p}=[p_n, n \in \mathcal{N}]$, $\boldsymbol{\theta}=[\theta_n, n \in \mathcal{N}]$, $\mathbf{f}=[f_{nm}, (n,m) \in \mathcal{L}]$, and $\mathbf{x}=[x_{nm}, (n,m) \in \mathcal{L_S}]$. With this notation in place, the DC-OTS problem can be formulated as follows:
\begin{subequations}\label{eq:OTS_NP}
\begin{align}
\min_{p_n,f_{nm},\theta_n,x_{nm}} & \quad \sum_{n} c_{n} \, p_{n} \label{eq:OTS_NP_obj}\\
\text{s.t.} & \quad \underline{p}_n \leq p_n \leq \overline{p}_n, \quad \forall n \in \mathcal{N} \label{eq:OTS_NP_Plimits}\\
&\quad \sum_{(n,m)\in\mathcal{L}^-_n}f_{nm} - \sum_{(n,m)\in\mathcal{L}^+_n}f_{nm} = p_n - d_n, \quad \forall n \in \mathcal{N} \label{eq:OTS_NP_PB}\\
& \quad f_{nm} = x_{nm}b_{nm}(\theta_n-\theta_m), \quad \forall (n,m) \in \mathcal{L}_\mathcal{S} \label{eq:OTS_NP_Flow_S}\\
& \quad f_{nm} = b_{nm}(\theta_n-\theta_m), \quad \forall (n,m) \in \mathcal{L} \setminus \mathcal{L}_\mathcal{S}\label{eq:OTS_NP_Flow_NS}\\
& \quad -x_{nm}\overline{f}_{nm}\leq f_{nm} \leq x_{nm}\overline{f}_{nm}, \quad \forall (n,m) \in \mathcal{L}_\mathcal{S} \label{eq:OTS_NP_Flow_limit_S}\\
& \quad -\overline{f}_{nm}\leq f_{nm} \leq \overline{f}_{nm}, \quad \forall (n,m) \in \mathcal{L} \setminus \mathcal{L}_\mathcal{S}\label{eq:OTS_NP_Flow_limit_NS}\\
& \quad x_{nm} \in \{0,1\}, \quad \forall (n,m) \in \mathcal{L}_\mathcal{S}\label{eq:OTS_NP_binary}\\
& \quad \theta_1 = 0 \label{eq:OTS_NP_slack}
\end{align}
\end{subequations}

The objective is to minimize the electricity generation cost, expressed as in~\eqref{eq:OTS_NP_obj}. For this, the power system operator essentially decides the lines that are switched off and the power output of generating units, which must lie within the interval $[\underline{p}_n,\overline{p}_n]$, as imposed in~\eqref{eq:OTS_NP_Plimits}. The flows through the transmission lines are governed by the so-called \emph{Kirchhoff's laws}, which translate into the nodal power balance equations~\eqref{eq:OTS_NP_PB} and the flow-angle relationship stated in~\eqref{eq:OTS_NP_Flow_S} and~\eqref{eq:OTS_NP_Flow_NS}. In the case of a switchable line, this relationship must be enforced only when the line is in service. This is why the binary variable $x_{nm}$ appears in~\eqref{eq:OTS_NP_Flow_S}. Naturally, $x_{nm} = 0$ must imply $f_{nm} = 0$. Constraints~\eqref{eq:OTS_NP_Flow_limit_S} and~\eqref{eq:OTS_NP_Flow_limit_NS} impose the capacity limits of the switchable and non-switchable lines, respectively. Constraint~\eqref{eq:OTS_NP_binary} states the binary character of variables $x_{nm}$, while equation~\eqref{eq:OTS_NP_slack} arbitrarily sets one of the nodal angles to zero to avoid solution multiplicity. 

Problem~\eqref{eq:OTS_NP} is a mixed-integer nonlinear programming problem due to the product $x_{nm}(\theta_n-\theta_m)$ in~\eqref{eq:OTS_NP_Flow_S}. This problem has been proven to be NP-hard even when the power network includes a spanning subnetwork connected by non-switchable lines only~\cite{fattahi2019bound} or takes the form of a series-parallel graph with a single generator/load pair~\cite{kocuk2016cycle}. The disjunctive nature of Equation~\eqref{eq:OTS_NP_Flow_S} allows for a linearization of Problem~\eqref{eq:OTS_NP} at the cost of introducing a pair of large enough constants $\underline{M}_{nm}$, $\overline{M}_{nm}$ per switchable line \cite{hedman2012flexible}. Indeed, Equation~\eqref{eq:OTS_NP_Flow_S} can be replaced by the inequalities
\begin{equation}\label{eq:linear}
\quad b_{nm}(\theta_n-\theta_m)-\overline{M}_{nm}(1-x_{nm}) \leq f_{nm} \leq b_{nm}(\theta_n-\theta_m)-\underline{M}_{nm}(1-x_{nm})
\end{equation}
provided that the large constants $\underline{M}_{nm}, \overline{M}_{nm}$ respectively constitute a lower and an upper bound of  $b_{nm}(\theta_n-\theta_m)$ when the line $(n,m)$ is disconnected ($x_{nm} = 0$), that is,
\begin{subequations}\label{eq:max}
\begin{align}
\underline{M}_{nm} &\leq \underline{M}^{\rm OPT}_{nm}:=b_{nm} \times\underset{\mathcal{F}}{\min} (\theta_n-\theta_m) \label{eq:max_min}\\
\overline{M}_{nm} &\geq \overline{M}^{\rm OPT}_{nm}:=b_{nm} \times\underset{\mathcal{F}}{\max} (\theta_n-\theta_m)\label{eq:max_max}
\end{align}
\end{subequations}
where $\mathcal{F}:= \{( \mathbf{p}, \boldsymbol{\theta}, \mathbf{f}, \mathbf{x}) \in \mathbb{R}^{2\vert\mathcal{N}\vert +\vert\mathcal{L}\vert + \vert\mathcal{L}_{\mathcal{S}}\vert}$ satisfying \eqref{eq:OTS_NP_Plimits}, \eqref{eq:OTS_NP_PB}, \eqref{eq:OTS_NP_Flow_NS}-\eqref{eq:OTS_NP_slack}, $x_{nm} = 0$, and $\eqref{eq:linear} \text{ for all } (n',m') \in \mathcal{L}_{\mathcal{S}}\setminus(n,m)$\}.  Note that, if $x_{nm} = 1$, we have
\[\quad b_{nm}(\theta_n-\theta_m) \leq f_{nm} \leq b_{nm}(\theta_n-\theta_m) \iff f_{nm} = b_{nm}(\theta_n-\theta_m).\]
Otherwise, i.e., if $x_{nm} = 0$, Equation~\eqref{eq:OTS_NP_Flow_limit_S} leads to $f_{nm} = 0$, which, together with~\eqref{eq:linear}, results in
\[\quad b_{nm}(\theta_n-\theta_m)-\overline{M}_{nm} \leq 0 \leq b_{nm}(\theta_n-\theta_m)-\underline{M}_{nm} \]
or, equivalently, 
\[\underline{M}_{nm}\le b_{nm}(\theta_n-\theta_m) \le \overline{M}_{nm}.\]
Finally, by Equation (3), we have
\[\underline{M}_{nm} \le \underline{M}^{\rm OPT}_{nm} \le b_{nm}(\theta_n-\theta_m) \le \overline{M}^{\rm OPT}_{nm} \le \overline{M}_{nm}.\]
First of all, for \eqref{eq:max} to be of any use, $\underline{M}^{\rm OPT}_{nm}$ and $\overline{M}^{\rm OPT}_{nm}$ must be finite. As proven in \cite{fattahi2019bound}, this is not the case in power systems where switching off lines can result in disconnected subnetworks. The possibility of islanding renders the minimization \eqref{eq:max_min} and the maximization \eqref{eq:max_max} unbounded. Consequently, the linearization of the DC-OTS problem based on \eqref{eq:linear} is not equivalent to its original nonlinear mixed-integer
formulation~\eqref{eq:OTS_NP} in this case. However, in practice, islanding in power grids is to be avoided in general for many reasons other than the minimization of the operational cost (e.g., due to reliability and security standards). Consequently, in what follows, we assume that the set of switchable lines $\mathcal{L}_{\mathcal{S}}$ is such that the connectivity of the whole power network is always guaranteed. In this setting, it is ensured that there exist finite valid large constants as stated in \eqref{eq:max}, namely, those corresponding to the longest path between nodes $n$ and $m$ on the undirected graph represented by the power grid. This already gives us an idea of how difficult the calculation of these constants is.
In this vein, the authors in~\cite{fattahi2019bound} show that, even when $\underline{M}^{\rm OPT}_{nm}$ and $\overline{M}^{\rm OPT}_{nm}$ are finite, computing them is as hard as solving the original DC-OTS problem. Therefore, we are obliged to be content with a lower and an upper bound. The choice of these bounds, or rather, of the large constants $\underline{M}_{nm}, \overline{M}_{nm}$ (for all $(n,m) \in \mathcal{L}_{\mathcal{S}}$) is of utmost importance, because it has a major impact on the relaxation bound of the mixed-integer \emph{linear} program that results from replacing~\eqref{eq:OTS_NP_Flow_S} with the inequalities~\eqref{eq:linear}, that is,
\begin{subequations} \label{subeq: MIP reformulation}
\begin{align}
\min_{p_n,f_{nm},\theta_n,x_{nm}} & \quad\sum_{n} c_{n} \, p_{n} \label{eq:ots_mip_of}\\
\text{s.t.} &\quad \eqref{eq:OTS_NP_Plimits}, \eqref{eq:OTS_NP_PB}, \eqref{eq:OTS_NP_Flow_NS}-\eqref{eq:OTS_NP_slack}\\
& \quad b_{nm}(\theta_n-\theta_m)-\overline{M}_{nm}(1-x_{nm}) \leq f_{nm}, \quad \forall (n,m) \in \mathcal{L}_{\mathcal{S}} \\
&\quad f_{nm} \leq b_{nm}(\theta_n-\theta_m)-\underline{M}_{nm}(1-x_{nm}), \quad \forall (n,m) \in \mathcal{L}_{\mathcal{S}} \\
& \quad f_{nm} = b_{nm}(\theta_n-\theta_m), \quad \forall (n,m) \in \mathcal{L} \setminus\mathcal{L}_{\mathcal{S}}
\end{align} \label{eq:ots_mip}
\end{subequations}

Tighter constants $\underline{M}_{nm}, \overline{M}_{nm}$ lead to stronger linear relaxations of~\eqref{subeq: MIP reformulation}, which, in turn, is expected to impact positively on the performance of the branch-and-cut algorithm used to solve it. Let us define $\mathbf{d}=[d_n, n \in \mathcal{N}]$ and $\mathbf{M}=[(\overline{M}_{nm}, \underline{M}_{nm}), (n,m) \in \mathcal{L_S}]$. We also define the lower and upper bounds of the binary decision variables as $\underline{\mathbf{x}} = [\underline{x}_{nm}, (n,m) \in \mathcal{L_S}]$ and $\overline{\mathbf{x}} = [\overline{x}_{nm}, (n,m) \in \mathcal{L_S}]$, respectively. Then, we denote as $\mathbf{x} = \text{OTS}(\mathbf{d},\mathbf{M}, \underline{\mathbf{x}},\overline{\mathbf{x}})$ the solution of model \eqref{eq:ots_mip} with the additional constraint $\underline{\mathbf{x}} \leq \mathbf{x} \leq \overline{\mathbf{x}}$. In the general case, $\underline{\mathbf{x}}=\mathbf{0}$ and $\overline{\mathbf{x}}=\mathbf{1}$. However, these bounds may change if the status of some switchable lines are fixed through learning.

On the assumption that the power network includes a spanning tree comprising non-switchable lines, the authors in~\cite{fattahi2019bound} propose the following symmetric bound: 
\begin{equation} \label{eq: big M Fattahi}
-\underline{M}_{nm} = \overline{M}_{nm} = b_{nm} \sum_{(k,l)\in {\rm SP}^0_{nm}}\frac{\overline{f}_{kl}}{b_{kl}}, \quad \forall (n,m) \in \mathcal{L}_{\mathcal{S}}
\end{equation}
where ${\rm SP}^0_{nm}$ is the shortest path between nodes $n$ and $m$ through said spanning tree. Note, however, that the shortest path between two nodes can be modified if some of the switchable lines are known to be connected. In that case, the resulting bounds are  reduced. Therefore, for a given status of the switchable lines $\mathbf{x}$, we denote by ${\rm SP}_{nm}(\mathbf{x})$ the updated shortest path, with ${\rm SP}^0_{nm} = {\rm SP}_{nm}(\mathbf{0})$. Besides, the bounds obtained using Equation \eqref{eq: big M Fattahi} with the updated shortest paths ${\rm SP}_{nm}(\mathbf{x})$ is referred to as $\mathbf{M} = \text{FAT}(\mathbf{x})$. This symmetric bound can be computed in polynomial time using Dijkstra's algorithm \cite{cormen2022introduction}.

In this paper, we propose and test simple, but effective data-driven scheme based on nearest neighbors to estimate lower bounds on $\underline{M}^{\rm OPT}_{nm}$ and upper bounds on $\overline{M}^{\rm OPT}_{nm}$. This scheme is also  used to fix some of the binaries $x_{nm}$ in \eqref{subeq: MIP reformulation}. While the inherent sampling error of the proposed methodology precludes optimality guarantees, our numerical experiments show that it is able to identify optimal or nearly-optimal solutions to the DC-OTS problem very fast. 

\section{Solution methods} \label{sec:met}

In this section, we present the different methods we consider to solve the DC-OTS problem. First, we describe the exact method proposed in \cite{fattahi2019bound}, which we use as a benchmark. Second, we explain a direct learning-based approach that utilizes the $K$ nearest neighbors technique and the learning-based heuristic approach investigated in \cite{hastie2009elements}. Finally, we introduce the data-based methodologies proposed in this paper. 

Suppose that the DC-OTS problem \eqref{eq:ots_mip} has been solved using the big-M values suggested in  \cite{fattahi2019bound} for different instances to form a training set $\mathcal{T} = \{(\mathbf{d}^t, \mathbf{x}^t, \boldsymbol{\theta}^t), t= 1, \ldots, \vert\mathcal{T}\vert\}$, where, the symbol $\vert\mathcal{T}\vert$ indicates the cardinal of set $\mathcal{T}$. For each instance, $t$, $\mathbf{d}^t = [d_n^t, {n\in\mathcal{N}]}$ denotes the vector of nodal loads, $\mathbf{x}^t=[x^t_{nm}, (n,m)\in\mathcal{L}_{\mathcal{S}}]$ is the vector of optimal binary variables, which determine whether line $(n, m)$ in instance $t$ is connected or not; and $\boldsymbol{\theta}^t = [\theta^t_n, n\in\mathcal{N}]$ is the vector of optimal voltage angles. 
For notation purposes, we use $C(\mathbf{d}^t, \mathbf{x}^t)$ to denote the value of the objective function \eqref{eq:OTS_NP_obj} when model \eqref{eq:OTS_NP} is solved for demand values $\mathbf{d}^t$ and the binary variables fixed to $\mathbf{x}^t$. This function can be evaluated for any set of feasible binary variables $\mathbf{x}^t$ by solving a linear programming problem. If this linear problem is infeasible, then $C(\mathbf{d}^t, \mathbf{x}^t) = \infty$. Additionally, for a given subset of instances $\mathcal{T}' \subset \mathcal{T}$, we define $\mathbf{x}(\mathcal{T}')$ as the component-wise average of the binary variables corresponding to the instances in $\mathcal{T}'$.

In what follows, we present different strategies to solve the DC-OTS problem for an unseen test instance $\hat{t}$ with demand values $\mathbf{d}^{\hat{t}}$. The goal is to employ the information from the training set, $\mathcal{T}$, to reduce the computational burden of solving the DC-OTS reformulation \eqref{subeq: MIP reformulation} for the test instance $\hat{t}$. Note that depending on the strategy that is applied, the response variable of the test instance to be learned can be $\mathbf{x}^{\hat{t}}$, $\boldsymbol{\theta}^{\hat{t}}$ or the tuple $(\mathbf{x}^{\hat{t}}, \boldsymbol{\theta}^{\hat{t}})$.

\subsection{Exact benchmark approach} \label{subs: benchmark method}

In the benchmark approach (\emph{Bench}) the optimal solution of the test DC-OTS problem is obtained using the proposal in \cite{fattahi2019bound}. Particularly, problem \eqref{subeq: MIP reformulation} is solved using the big-M values computed according to Equation \eqref{eq: big M Fattahi}. This strategy is an exact approach that does not make use of previously solved instances of the problem, but guarantees that its global optimal solution is eventually retrieved. Nevertheless, the computational time employed by this approach may be extremely high. Algorithm \ref{alg:bench} shows a detailed description of this approach.

\begin{algorithm}
\caption{\emph{Bench}}\label{alg:bench}
\textbf{Input:} load vector for test instance $\hat{t}$, $\mathbf{d}^{\hat{t}}$.\\
\begin{enumerate}[label={\arabic*)}]
\item Compute the bounds with all switchable lines open, i.e., $\mathbf{M}^0 = \text{FAT}(\mathbf{0})$.\\
\item Solve $\mathbf{x}^{\hat{t}} = \text{OTS}(\mathbf{d}^{\hat{t}},\mathbf{M}^0,\mathbf{0},\mathbf{1})$ \\
\end{enumerate}
\textbf{Output:} Optimal network configuration $\mathbf{x}^{\hat{t}}$.
\end{algorithm}

\subsection{Existing learning-based approaches} \label{subs:existing}

In this subsection we present two existing learning approaches based on the $K$ nearest neighbors technique \cite{hastie2009elements}. The first approach is a pure machine-learning strategy that directly predicts the binary variables of the test instance using the information of the $K$ closest training data. Such closeness is measured in terms of the $\ell_2$ distance among the load values of the training and test points, that is, by computing $\|\mathbf{d}^t-\mathbf{d}^{\hat{t}}\|_2$, for $t = 1, \ldots, \vert\mathcal{T}\vert\}$. For each test instance $\hat{t}$, the set of $K$ closest instances is denoted as  $\mathcal{T}^{\hat{t}}_K = K\text{NN}(\mathbf{d}^{\hat{t}})$. This method is referred to as \emph{Direct} since it \textit{directly} predicts the value of all binary variables from the data. 

In the particular case of the DC-OTS problem, we adapt the $K$nn strategy as follows: for a fixed number of neighbors $K$, we fix the binary variables of the test problem \eqref{eq:OTS_NP} to the rounded mean of the decision binary variables of such $K$ nearest neighbors. Once all binary variables are fixed, model \eqref{eq:OTS_NP} becomes a linear programming problem that can be rapidly solved. Algorithm \ref{alg:direct} shows a detailed explanation of the procedure. Note that, in this strategy we only need the information about the load vector and the optimal binary variables in the training data, i.e., we only need $\{(\mathbf{d}^t, \mathbf{x}^t)\}$ for $t=1, \ldots, \vert\mathcal{T}\vert\}$. 
This approach is very simple and fast. However, fixing the binary variables using a rounding procedure may yield a non-negligible number of infeasible and suboptimal problems.

\begin{algorithm}
\caption{\emph{Direct}}\label{alg:direct}
\textbf{Input:} number of neighbors, $K$; training set, $\mathcal{T} = \{(\mathbf{d}^t, \mathbf{x}^t)\}$ for $t= 1, \ldots, \vert\mathcal{T}\vert\}$; and load vector for test instance $\hat{t}$, $\mathbf{d}^{\hat{t}}$.\\
\begin{enumerate}[label={\arabic*)}]
\item $\mathcal{T}^{\hat{t}}_K = K\text{NN}(\mathbf{d}^{\hat{t}})$. \\
\item \label{item: step 3 algorithm $K$nn} Compute the binary variables $\mathbf{x}^{\hat{t}} = \lceil\mathbf{x}(\mathcal{T}^{\hat{t}}_K)\rfloor$, where $\lceil\mathbf{x}\rfloor$ denotes the component-wise nearest integer function.\\
\end{enumerate}
\textbf{Output:} Network configuration $\mathbf{x}^{\hat{t}}$.
\end{algorithm}

The second learning-based methodology explained in this subsection is proposed in \cite{johnson2020knearest} and also employs the $K$nn technique. As occurs in the previous strategy, here, the authors assume given the set $\{(\mathbf{d}^t, \mathbf{x}^t),\, t = 1, \ldots, \vert\mathcal{T}\vert\}$.
In short, their proposal works as follows: for a fixed value of $K$, the $K$ closest instances to the test point are saved in the set $\mathcal{T}^{\hat{t}}_K$. Then, we evaluate function $C(\mathbf{d}^{\hat{t}},\mathbf{x}^t)$ for each $t\in\mathcal{T}^{\hat{t}}_K$ by solving $K$ linear problems. The optimal binary variables for the test instance $\mathbf{x}^{\hat{t}}$ are set to those $\mathbf{x}^t$ that lead to the lowest value of $C(\mathbf{d}^{\hat{t}},\mathbf{x}^t)$. This approach is denoted as \emph{Linear} and more details about it are provided in Algorithm \ref{alg:linear}.

\begin{algorithm}
\caption{\emph{Linear}  \cite{johnson2020knearest}}\label{alg:linear}
\textbf{Input:} number of neighbors, $K$; training set, $\mathcal{T} = \{(\mathbf{d}^t, \mathbf{x}^t), \, t= 1, \ldots, \vert\mathcal{T}\vert\}$; and load vector for test instance $\hat{t}$, $\mathbf{d}^{\hat{t}}$.\\
\begin{enumerate}[label={\arabic*)}]
\item $\mathcal{T}^{\hat{t}}_K = K\text{NN}(\mathbf{d}^{\hat{t}})$. \\
\item Select $\tilde{t} = \arg\min\limits_{t\in\mathcal{T}^{\hat{t}}_K} C(\mathbf{d}^{\hat{t}},\mathbf{x}^t)$. \\
\end{enumerate}
\textbf{Output:} Network configuration $\mathbf{x}^{\tilde{t}}$.
\end{algorithm}

Note that the value of $K$ strongly affects the speed of the algorithm as well as the number of suboptimal or infeasible problems. Larger values of $K$ imply taking into account more training points to get the estimation of the test response. As a consequence, a larger number of LPs should be solved, and the computational burden increases. However, the probability of having suboptimal or, even worse, infeasible solutions is reduced. On the contrary, lower values of $K$ diminishes the computational time of the procedure but increases the risk of obtaining suboptimal or infeasible solutions.

\subsection{Proposed learning-based approaches} \label{subs:proposed}

In this subsection, we propose two improved methodologies which combine the benefits of exact and learning methods. Both approaches start by finding the $K$ closest training points to the test instance $\hat{t}$ and fixing those binary variables that reach the same value for all nearest neighbors according to an \textit{unanimous vote}. The two proposed approaches also find, in a different fashion, lower values of the big-Ms than those computed in \cite{fattahi2019bound}. Since some binary variables may have been fixed to one thanks to the neighbors' information, the first approach we propose consists in recomputing the shortest paths and the corresponding big-M values using \eqref{eq: big M Fattahi}. Differently, the second methodology proposed in this paper directly set the big-M values to the maximum and minimum values of the angle differences observed in the closest DC-OTS instances. Either way, smaller big-Ms are obtained, and hence, the associated feasible region of the DC-OTS problem is tighter. As a consequence, we solve a single MILP with a tighter feasible region and a smaller number of binary variables.

More specifically, in the first proposed approach  (denoted as \emph{FixB-FatM}) the binary variables of the test instance are set to 1 (resp. to 0) if all the training instances in $\mathcal{T}^{\hat{t}}_K$ concur that the value should be 1 (resp. 0). On the other hand, for those binary variables that are not fixed, the corresponding big-M values are updated using the information of the previously fixed variables. In particular, these fixed binaries are used to recompute the shortest path that determines the big-M values in Equation \eqref{eq: big M Fattahi}. In essence, the computation of the new shortest path involves not only the non-switchable lines from the original spanning tree but also those switchable lines with a learned status equal to 1. This update could result in even shorter paths, leading to improved big-M bounds and a more tightly defined feasible region. This strategy relies on the unanimity of all the nearest neighbors and therefore, this learning-based approach is expected to be quite conservative, specially for high values of $K$.  

In order to further assess the computational savings yielded by this approach we also investigate two variations. For instance, we denote by \emph{FixB} the approach in which binary variables are fixed but big-M values are computed using only the information from the original spanning tree. We also consider the \emph{FatM} approach that does not fix any binary decision variable but only uses the information of the closest neighbors to recompute the shortest paths and update the big-M values with Equation \eqref{eq: big M Fattahi}. In other words, while none of the binary variables are fixed in this method, the learned status of switchable lines can still be utilized to decrease the big-M values. By comparing the computational burden of these three approaches we can analyze whether the numerical improvements are caused by the lower number of binary variables or the tighter values of the big-M parameters. Algorithms \ref{alg:fixb}, \ref{alg:fatm} and \ref{alg:fitb-fatm} show a detailed description of the methods \emph{FixB}, \emph{FatM} and \emph{FixB-FatM}, respectively. 

\begin{algorithm}
\caption{\emph{FixB}}\label{alg:fixb}
\textbf{Input:} number of neighbors, $K$; training set, $\mathcal{T} = \{(\mathbf{d}^t, \mathbf{x}^t)\}$ for $t= 1, \ldots, \vert\mathcal{T}\vert$; and load vector for test instance $\hat{t}$, $\mathbf{d}^{\hat{t}}$.\\
\begin{enumerate}[label={\arabic*)}]
\item $\mathbf{M}^0 = \text{FAT}(\mathbf{0})$ and  $\mathcal{T}^{\hat{t}}_K = K\text{NN}(\mathbf{d}^{\hat{t}})$. \\
\item Compute $\underline{\mathbf{x}}^{\hat{t}} = \lfloor\mathbf{x}(\mathcal{T}^{\hat{t}}_K)\rfloor$ and $\overline{\mathbf{x}}^{\hat{t}} = \lceil\mathbf{x}(\mathcal{T}^{\hat{t}}_K)\rceil$.\\
\item Solve $\mathbf{x}^{\hat{t}} = \text{OTS}(\mathbf{d}^{\hat{t}},\mathbf{M}^0,\underline{\mathbf{x}}^{\hat{t}},\overline{\mathbf{x}}^{\hat{t}})$.\\
\end{enumerate}
\textbf{Output:} Network configuration $\mathbf{x}^{\hat{t}}$.
\end{algorithm}

\begin{algorithm}
\caption{\emph{FatM}}\label{alg:fatm}
\textbf{Input:} number of neighbors, $K$; training set, $\mathcal{T} = \{(\mathbf{d}^t, \mathbf{x}^t)\}$ for $t= 1, \ldots, \vert\mathcal{T}\vert$; and load vector for test instance $\hat{t}$, $\mathbf{d}^{\hat{t}}$.\\
\begin{enumerate}[label={\arabic*)}]
\item $\mathcal{T}^{\hat{t}}_K = K\text{NN}(\mathbf{d}^{\hat{t}})$. \\
\item Compute $\widetilde{\mathbf{x}}^{\hat{t}} = \lfloor\mathbf{x}(\mathcal{T}^{\hat{t}}_K)\rfloor$ and $\widetilde{\mathbf{M}}^{\hat{t}} = \text{FAT}(\widetilde{\mathbf{x}}^{\hat{t}})$.\\
\item Solve $\mathbf{x}^{\hat{t}} = \text{OTS}(\mathbf{d}^{\hat{t}},\widetilde{\mathbf{M}}^{\hat{t}},\mathbf{0},\mathbf{1})$.\\
\end{enumerate}
\textbf{Output:} Network configuration $\mathbf{x}^{\hat{t}}$.
\end{algorithm}

\begin{algorithm}
\caption{\emph{FixB-FatM}}\label{alg:fitb-fatm}
\textbf{Input:} number of neighbors, $K$; training set, $\mathcal{T} = \{(\mathbf{d}^t, \mathbf{x}^t)\}$ for $t= 1, \ldots, \vert\mathcal{T}\vert$; and load vector for test instance $\hat{t}$, $\mathbf{d}^{\hat{t}}$.\\
\begin{enumerate}[label={\arabic*)}]
\item $\mathcal{T}^{\hat{t}}_K = K\text{NN}(\mathbf{d}^{\hat{t}})$. \\
\item Compute $\widetilde{\mathbf{x}}^{\hat{t}} = \lfloor\mathbf{x}(\mathcal{T}^{\hat{t}}_K)\rfloor$ and $\widetilde{\mathbf{M}}^{\hat{t}} = \text{FAT}(\widetilde{\mathbf{x}}^{\hat{t}})$.\\
\item Determine $\underline{\mathbf{x}}^{\hat{t}} = \lfloor\mathbf{x}(\mathcal{T}^{\hat{t}}_K)\rfloor$ and $\overline{\mathbf{x}}^{\hat{t}} = \lceil\mathbf{x}(\mathcal{T}^{\hat{t}}_K)\rceil$.\\
\item Solve $\mathbf{x}^{\hat{t}} = \text{OTS}(\mathbf{d}^{\hat{t}},\widetilde{\mathbf{M}}^{\hat{t}},\underline{\mathbf{x}}^{\hat{t}},\overline{\mathbf{x}}^{\hat{t}})$.\\
\end{enumerate}
\textbf{Output:} Network configuration $\mathbf{x}^{\hat{t}}$.
\end{algorithm}

\textit{FixB} and \textit{FixB-FatM} can be slightly modified to relax the unanimity condition required to fix binary variables. To do so, we introduce a threshold parameter $\tau < 0.5$. The binary variables are then fixed according to the following rules:
\begin{itemize}
\item[-] If the predicted status for a particular line falls in $[0,\tau]$, the binary variable is fixed to 0.
\item[-] If the predicted status for a particular line falls in $[1-\tau, 1]$, the binary variable is fixed to 1.
\item[-] If the predicted status for a particular line falls between $(\tau, 1-\tau)$, the binary variable is left unfixed.
\end{itemize}
This can be implemented by replacing, respectively, step 2) in Algorithm \ref{alg:fixb} and step 3) in Algorithm \ref{alg:fitb-fatm} by:
\[\underline{\mathbf{x}}^{\hat{t}} = \lfloor\min(\mathbf{x}(\mathcal{T}^{\hat{t}}_K)+\tau,1)\rfloor \text{ and } \overline{\mathbf{x}}^{\hat{t}} = \lceil\max(\mathbf{x}(\mathcal{T}^{\hat{t}}_K)-\tau,0)\rceil\]

The value of $K$ also plays an important role in these approaches. Low values of $K$ increase the chances of unanimous consensus of the nearest neighbors and therefore, a higher number of binary variables are expected to be fixed, and tighter big-M values are obtained. This way, the computational burden of the OTS problem is reduced at the expense of increasing the risk of obtaining infeasible or suboptimal problems. In the extreme case, if $K=1$, all binary variables are fixed to the values of the closest instance of the training set. On the contrary, large values of $K$ increase the computational burden but the resulting problems have a high chance of being feasible. In the extreme case, if the whole training set is considered, very few binary variables are expected to be fixed and the computational savings are reduced.

The three methodologies presented above compute the big-M values using past observed data through the shortest path algorithm. However, as can be derived from Equation \eqref{eq:max}, the values $\overline{M}_{nm}$ and $\underline{M}_{nm}$ for a switchable line are just the maximum and minimum values of the difference between the voltage angles at nodes $n$ and $m$ multiplied by $b_{nm}$. Therefore, following this idea, the second data-driven approach that we propose (denoted as \emph{FixB-AngM}) estimate the big-M values using the information of historic observed angles as follows:
\begin{equation} \label{eq:bigm_ang}
\overline{M}_{nm}= b_{nm}\times 
\underset{t\in \mathcal{T}:\, x_{nm}^t= 0}{\max} (\theta_n^t - \theta_m^t) \qquad \underline{M}_{nm}= b_{nm}\times 
\underset{t\in \mathcal{T}:\, x_{nm}^t= 0}{\min}  (\theta_n^t - \theta_m^t)
\end{equation}
Using \eqref{eq:bigm_ang} to compute the bounds values for a set of past instances $\mathcal{T}$ is denoted as $\mathbf{M}=\text{ANG}(\mathcal{T})$ for notation purposes. It is important to clarify that computing the big-M values using \eqref{eq:max} and \eqref{eq:bigm_ang} involves significant differences. The problems addressed by \eqref{eq:max} focus on identifying the tightest valid bounds by solving mixed-integer problems, which are as challenging as the original OTS problem. In contrast, Equation \eqref{eq:bigm_ang} efficiently approximates these bounds using observed angles from the historical dataset. Consequently, the bounds derived from \eqref{eq:bigm_ang} are consistently tighter than those obtained from \eqref{eq:max}, potentially excluding feasible solutions to the original OTS problem if the training set lacks sufficient representativeness. In fact, this strategy is riskier than the one used in \emph{FixB-FatM} since it leads to much tighter feasible regions, which significantly reduces the computational burden of solving the OTS problem, but also increases the chances of yielding infeasible problems. To avoid using too tight big-M values that could cut off the optimal solution, the learned bounds obtained through \eqref{eq:bigm_ang} can be multiplied by a security factor $\lambda\geq 1$.

For the sake of comparison, we also consider the approach \emph{AngM} in which no binary variables are fixed and big-M values are set using the observed angle differences. More details about the  approaches \emph{FixB-AngM} and \emph{AngM} are provided in Algorithms \ref{alg:fixb_angm} and \ref{alg:angm}, respectively. It is worth noticing that while the big-M values computed by \eqref{eq: big M Fattahi} are symmetric, those derived by Algorithms \ref{alg:fixb_angm} and \ref{alg:angm} are not. 

\begin{algorithm}
\caption{\emph{FixB-AngM}}\label{alg:fixb_angm}
\textbf{Input:} number of neighbors, $K$; training set, $\mathcal{T} = \{(\mathbf{d}^t, \mathbf{x}^t, \boldsymbol{\theta}^t)\}$ for $t= 1, \ldots, \vert\mathcal{T}\vert$; load vector for test instance $\hat{t}$, $\mathbf{d}^{\hat{t}}$; and security factor $\lambda \geq 1$.\\
\begin{enumerate}[label={\arabic*)}]
\item $\mathcal{T}^{\hat{t}}_K = K\text{NN}(\mathbf{d}^{\hat{t}})$. \\
\item Determine $\underline{\mathbf{x}}^{\hat{t}} = \lfloor\mathbf{x}(\mathcal{T}^{\hat{t}}_K)\rfloor$ and $\overline{\mathbf{x}}^{\hat{t}} = \lceil\mathbf{x}(\mathcal{T}^{\hat{t}}_K)\rceil$.\\
\item Compute big-M values as $\widehat{\mathbf{M}}=  \lambda\cdot \text{ANG}(\mathcal{T})$.\\
\item Solve $\mathbf{x}^{\hat{t}} = \text{OTS}(\mathbf{d}^{\hat{t}},\widehat{\mathbf{M}},\underline{\mathbf{x}}^{\hat{t}},\overline{\mathbf{x}}^{\hat{t}})$.\\
\end{enumerate}
\textbf{Output:} Network configuration $\mathbf{x}^{\hat{t}}$.
\end{algorithm}

\begin{algorithm}
\caption{\emph{AngM}}\label{alg:angm}
\textbf{Input:} training set, $\mathcal{T} = \{(\mathbf{d}^t, \mathbf{x}^t, \boldsymbol{\theta}^t)\}$ for $t= 1, \ldots, \vert\mathcal{T}\vert$; load vector for test instance $\hat{t}$, $\mathbf{d}^{\hat{t}}$; and security factor $\lambda \geq 1$.\\
\begin{enumerate}[label={\arabic*)}]
\item Compute big-M values as $\widehat{\mathbf{M}}= \lambda\cdot\text{ANG}(\mathcal{T})$.\\
\item Solve $\mathbf{x}^{\hat{t}} = \text{OTS}(\mathbf{d}^{\hat{t}},\widehat{\mathbf{M}},\mathbf{0},\mathbf{1})$.\\
\end{enumerate}

\textbf{Output:} Network configuration $\mathbf{x}^{\hat{t}}$.
\end{algorithm}

To sum up, Table \ref{tab: summary methods} provides a brief description of the different methods explained throughout Section \ref{sec:met}. The first column of the table includes the name of each strategy.  The second column shows whether the final problem to be solved is a linear program (LP) or a mixed-integer linear program (MILP). In the third column, the total number of problems to be solved is indicated. Column four shows the number of binary decision variables of the MILPs to be solved. Particularly, \emph{original} means that the number of variables is exactly the same as the one from the original OTS formulation \eqref{subeq: MIP reformulation}. In contrast, \emph{reduced} means that the number of binary variables of the resulting MILP has been reduced compared to the original formulation. Finally, the last column indicates how the big-M values have been computed. If \emph{shortest (spanning)} is written, then we indicate that the bounds are computed by means of the shortest path method and only using the information from the original spanning subgraph. On the contrary, the choice \emph{shortest (update)} means that the shortest paths needed to compute the big-M values have been updated with the information provided by the closest neighbors. Finally, the word \emph{historic angles} implies that the bounds are computed using the voltage angle information of previously solved instances.

\begin{table}[ht]
\centering
\begin{tabular}{lcccccc}
\toprule
Method  & LP/MILP & \# problems & \# binary & big-M computation\\
\midrule
 \emph{Bench}  & MILP &1& original & shortest (spanning) \\
 \emph{Direct} &-& - & - & - \\
 \emph{Linear}& LP &K& -  & - \\
 \emph{FixB}  & MILP &1& reduced & shortest (spanning) \\
 \emph{FatM}  & MILP &1& original & shortest (update) \\
 \emph{FixB-FatM} & MILP &1& reduced & shortest (update) \\
 \emph{FixB-AngM}  & MILP &1& reduced & historic angles\\
 \emph{AngM}  & MILP &1& original & historic angles\\
\bottomrule
\end{tabular}
\caption{Summary of the methods explained in Section \ref{sec:met}.}
\label{tab: summary methods}
\end{table}

\section{Numerical simulations} \label{sec: num}

In this section, we present the computational results of the different methodologies discussed in Section \ref{sec:met} for a realistic network. In particular, we compare all approaches using a 118-bus network that includes 186 lines \cite{blumsack2006network}. This network size is sufficiently substantial to render the instances nontrivial for current algorithms, yet not so large as to make them computationally intractable. Indeed, this is the most commonly used network to test OTS solving strategies in the literature \cite{fisher2008optimal, kocuk2016cycle, fattahi2019bound, johnson2020knearest, dey2022node}. As justified in Section \ref{sec:ots}, we consider a fixed connected spanning subgraph of 117 lines, while the remaining 69 lines can be switched on or off to minimize the operation cost. The spanning subgraph has been chosen in order to obtain sufficiently challenging problems. For this network, we generate 500 different instances of the OTS problem that differ in the nodal demand $d_n$ using probability distributions centered in the baseline demand $\widehat{d}_n$. Since the demand variability may significantly affect the performance of the compared methodologies, we consider the following three cases:
\begin{itemize}
    \item[-] \emph{Unif10}: The demand levels are sampled using independent uniform distributions in the range $[0.9\widehat{d}_n,1.1\widehat{d}_n]$.
    \item[-] \emph{Unif20}: The demand levels are sampled using independent uniform distributions in the range $[0.8\widehat{d}_n,1.2\widehat{d}_n]$.
    \item[-] \emph{Normal}: The demand levels are sampled using a multinormal distribution with the correlation matrix obtained from the demand time series available at \cite{joswig2021opflearndata}.
\end{itemize}
The three database files can be downloaded from \cite{OASYS2023learning}. We use a leave-one-out cross-validation technique under which all the available data except for one data point is used as the training set, and the left-out data point is used as the test set. Consequently, the number of nearest neighbors $K$ ranges from 1 to 499. This process is repeated for all data points and the resulting performance metrics are averaged to get an estimate of the model's generalization performance.

All optimization problems have been solved using GUROBI 9.1.2 \cite{gurobi} on a Linux-based server with CPUs clocking at 2.6 GHz, 1 thread and 8 GB of RAM. In all cases, the optimality GAP has been set to $0.01\%$ and the time limit to 1 hour. 

The simulation results are presented in two subsections. In Subsection \ref{sub:base_case_study}, a comprehensive comparison is conducted for all learning strategies introduced in Section \ref{sec:met} using the \emph{Unif10} database. Subsection \ref{sub:impact_variability} utilizes the \emph{Unif20} and \emph{Normal} databases to explore the impact of increased demand variability and correlation on the computational performance of these methodologies. 

\subsection{Base case study}\label{sub:base_case_study}

All simulation results presented in this subsection correspond to the \emph{Unif10} database. To illustrate the economic advantages of disconnecting some lines, Figure~\ref{fig:benchmark_savings} depicts an histogram of the relative difference between the DC-OTS cost if model \eqref{eq:ots_mip} is solved by the benchmark approach described in Section \ref{subs: benchmark method}, and the cost obtained if all the 186 lines are connected. This second cost is computed by fixing binary variables $x_{nm}$ to one and solving model \eqref{eq:OTS_NP} as a linear programming problem. Figure~\ref{fig:benchmark_savings} does not include the instances for which this linear problem is infeasible. As observed, the cost savings are significant in most instances, and in the most favorable cases it reaches 15\%. The average cost savings for this particular network and the 500 instances considered is 13.2\%. On the other hand, solving model \eqref{eq:ots_mip} is computationally hard and to prove it, Figure~\ref{fig:benchmark_time} plots the number of problems solved as a function of the computational time.  For illustration purposes, the left plot shows the 439 problems solved in less than 100 seconds (``easy'' instances) and the right plot the remaining ``hard'' instances that require a longer time. The average time of all instances is 145s, while the average time of the hard instances amounts to 1085 seconds, which demonstrates the difficulty of solving model \eqref{eq:ots_mip} to certified optimality. In addition, the benchmark approach is unable to solve 12 of the 500 instances to global optimality within one hour (with a maximum mip-gap equal to 2.46\%) even though model \eqref{eq:ots_mip} ``only'' includes 69 binary variables associated to the 69 switchable lines. This means that, for these 12 instances, this method has not been able to certify the optimality of the best integer solution found within the time limit, due to the poor relaxation bound originated from excessively large big-M values. We have thoroughly examined the simulations of this case study and verified that, for all instances, the best integer solution identified by the benchmark consistently matches the best solution discovered by all the other (learning-based) approaches. This lead us to conjecture that the benchmark does find the optimal solution for all instances in the \emph{Unif10} database. Therefore, throughout this section, we compare the different methodologies with the best integer solution found in one hour by the \emph{Bench} approach.

\begin{figure}[ht]
\centering
\begin{tikzpicture}[scale=0.65]
\begin{axis}[ymin=0, ymax=50, xmin=10, xmax=15.2, minor y tick num = 3, area style,xlabel=Cost savings (\%),ylabel=\# instances]
\addplot+[ybar interval,mark=no] plot coordinates {(10.04,1.0)(10.21,0.0)(10.38,0.0)(10.55,0.0)(10.72,0.0)(10.89,1.0)(11.05,2.0)(11.22,2.0)(11.39,3.0)(11.56,4.0)(11.73,8.0)(11.9,9.0)(12.07,13.0)(12.23,21.0)(12.4,21.0)(12.57,25.0)(12.74,30.0)(12.91,47.0)(13.08,29.0)(13.25,37.0)(13.41,32.0)(13.58,27.0)(13.75,21.0)(13.92,21.0)(14.09,16.0)(14.26,21.0)(14.43,8.0)(14.59,6.0)(14.76,0.0)(14.93,3.0)(15.1,0)};
\end{axis}
\end{tikzpicture}
\caption{DC-OTS cost savings distribution}
\label{fig:benchmark_savings}
\end{figure}

\begin{figure}[ht]
\captionsetup[subfigure]{labelformat=empty}
\centering
\begin{subfigure}[b]{0.45\textwidth}
\centering
\begin{tikzpicture}[scale=0.65]
\begin{axis}[ xlabel=Time (seconds), ylabel= \# problems solved,ymin=0, ymax=500, xmin=0, xmax=100, legend pos=south east, legend cell align={left}, extra x ticks = {100}] 
\addplot[color=black,line width=1] table[y=count, x=BEN, col sep=comma]{times.tex}; 
\end{axis} 
\end{tikzpicture}
\caption{}
\label{fig:benchmark_easy}
\end{subfigure}
\hfill
\begin{subfigure}[b]{0.45\textwidth}
\centering
\begin{tikzpicture}[scale=0.65]
\begin{axis}[ xlabel=Time (seconds), ylabel= \# problems solved,ymin=440, ymax=500, xmin=100, xmax=3600, legend pos=south east, legend cell align={left}, extra x ticks = {100}] 
\addplot[color=black,line width=1] table[y=count, x=BEN, col sep=comma]{times.tex}; 
\end{axis} 
\end{tikzpicture}
\caption{}
\label{fig:benchmark_hard}
\end{subfigure} 
\caption{Computational burden of the \emph{Bench} approach}
\label{fig:benchmark_time}
\end{figure}

Next, we discuss the results provided by the \emph{Direct} approach described in Section \ref{subs:existing}, where the binary variables are just fixed to the values  predicted by the nearest neighbor technique. Table \ref{tab:$K$nn_results} collates, for different number of neighbors $K$, the number of instances in which \emph{Direct} delivers the same solution obtained by the benchmark (\# opt), the number of instances with a suboptimal solution (\# sub) as well as the average and maximum relative gap with respect to the benchmark approach (gap-ave, gap-max). Note that the metrics \# opt,  \# sub, gap-ave and gap-max are computed with respect to the best solution found within one hour, which may not correspond to the true optimum. Finally, the results are also compared in terms of the average computational time, which can be seen in the last column of the table. Unsurprisingly, this approach is extremely fast and the computational time is just negligible. On the other hand, the vast majority of the instances only attain suboptimal solutions for any number of neighbors $K$, and the maximum gap is above 8\% in all cases. These results illustrate that the use of machine-learning approaches to directly predict the value of the binary variables of mixed-integer problems is likely to be extremely fast but potentially  suboptimal. 

\begin{table}[ht]
\setlength{\tabcolsep}{3pt}
\centering
\begin{tabular}{ccccccc}
\toprule
$K$ & \# opt & \# sub & gap-ave & gap-max & time (s) \\
\midrule
 5 & 2 & 498 & 1.799 & 13.78 & 0.0 \\
 10 & 0 & 500 & 2.025 & 16.40 & 0.0 \\
 20 & 0 & 500 & 2.085 & 13.84 & 0.0 \\
 50 & 0 & 500 & 2.057 & 14.13 & 0.0 \\
 100 & 0 & 500 & 1.846 & 12.53 & 0.0 \\
 200 & 0 & 500 & 2.367 & 12.28 & 0.0 \\
 499 & 0 & 500 & 2.629 & 8.38 & 0.0 \\
\bottomrule
\end{tabular}
\caption{Performance of the \emph{Direct} approach}
\label{tab:$K$nn_results}
\end{table}

Now we run similar experiments using the \emph{Linear} approach described in Section \ref{subs:existing} and proposed in \cite{johnson2020knearest}. The corresponding results are presented in Table \ref{tab:$K$nn_lp_results}. Logically, the \emph{Linear} solves a higher number of LP problems for different combinations of the binary variables and therefore, some instances are solved to optimality, specially for large values of $K$. Although this methodology could be parallelized, Table \ref{tab:$K$nn_lp_results} includes the sum of the computational times required to solve all the LP problems and therefore, this time increases with $K$. It is worth clarifying that the computational time required to find the nearest neighbors is below 1ms in all cases. Although the computational burden is insignificant if compared with the benchmark, the number of suboptimal cases and maximum gap are still considerable.

\begin{table}[ht]
\setlength{\tabcolsep}{3pt}
\centering
\begin{tabular}{ccccccc}
\toprule
 $K$ & \# opt & \# sub & gap-ave & gap-max & time (s) \\
\midrule
 5 & 10 & 490 & 0.300 & 3.59 & 0.00 \\
 10 & 19 & 481 & 0.194 & 3.56 & 0.01 \\
 20 & 24 & 476 & 0.130 & 1.61 & 0.02 \\
 50 & 51 & 449 & 0.083 & 1.06 & 0.04 \\
 100 & 77 & 423 & 0.061 & 0.71 & 0.08 \\
 200 & 104 & 396 & 0.049 & 0.71 & 0.16 \\
 499 & 127 & 373 & 0.043 & 0.71 & 0.39 \\
\bottomrule
\end{tabular}
\caption{Performance of the \emph{Linear} approach}
\label{tab:$K$nn_lp_results}
\end{table}

We continue this numerical study by comparing approaches \emph{FixB}, \emph{FatM} and \emph{FixB-FatM} discussed in Section \ref{subs:proposed}. For simplicity, Table \ref{tab:b_m_bm_results} provides the results for $K=50$ (10\% of the training data) and $\tau=0$. Unlike \emph{Direct} and \emph{Linear}, these three approaches lead to the optimal solution for all instances, which confirms their robustness for a sufficiently high number of neighbors. Therefore, although these approaches require a higher computational burden than \emph{Direct} and \emph{Linear}, they still involve significant computational savings with respect to the benchmark, while reducing the probability of returning suboptimal solutions. 

Table \ref{tab:b_m_bm_results} also shows that approaches \emph{FixB}, \emph{FatM} and \emph{FixB-FatM} differ in terms of their computational burden. The \emph{FatM} approach reports higher times than \emph{FixB}, which allows us to conclude that fixing some binary variables involves higher computational savings than tightening the big-M constants. Notwithstanding this, the highest computational gains are obtained if both effects are combined under the \emph{FixB-FatM} approach. Figure \ref{fig:b_m_bm_results} plots the number of problems solved as a function of time. In the left subplot, the x-axis ranges from 0 to 100s, while in the right subplot the x-axis goes from 100s to 3600s. In the left subplot we can observe that approaches \emph{FixB} and \emph{FixB-FatM} are able to solve most of the instances in less than 100s, while approach \emph{FatM} has a similar performance as the benchmark. In the right subplot we see that the hardest instance solved by \emph{FixB} and \emph{FixB-FatM} requires 1645s and 296s, respectively. On the contrary, although \emph{FatM} outperforms the benchmark, this approach is not able to solve all instances in less than one hour. 

\begin{table}[ht]
\setlength{\tabcolsep}{3pt}
\centering
\begin{tabular}{lccccc}
\toprule
& \# opt & \# sub & gap-ave & gap-max & time (s) \\
\midrule
 \emph{FixB} & 500 & 0 & - & - & 16.39 \\
 \emph{FatM} & 500 & 0 & - & - & 109.95 \\
 \emph{FixB-FatM} & 500 & 0 & - & - & 12.33 \\ 
 \bottomrule
\end{tabular}
\caption{Performance of \emph{FixB}, \emph{FatM}, \emph{FixB-FatM} for $K$=50 and $\tau = 0$}
\label{tab:b_m_bm_results}
\end{table} 

\begin{figure}[ht]
\captionsetup[subfigure]{labelformat=empty}
\centering
\begin{subfigure}[b]{0.45\textwidth}
\centering
\begin{tikzpicture}[scale=0.65]
\begin{axis}[ xlabel=Time (seconds), ylabel= \# problems solved,ymin=0, ymax=500, xmin=0, xmax=100, legend pos=south east, legend cell align={left}] 
\addplot[color=blue,line width=1] table[y=count, x=MIP_BM50, col sep=comma]{times.tex};
\addplot[color=green,line width=1,dashdotted] table[y=count, x=MIP_B50, col sep=comma]{times.tex}; 
\addplot[color=purple,line width=1,dashed] table[y=count, x=MIP_M50, col sep=comma]{times.tex}; 
\addplot[color=red,line width=1,dotted] table[y=count, x=BEN, col sep=comma]{times.tex}; 
\legend{\emph{FixB-FatM},\emph{FixB},\emph{FatM},\emph{Bench}}
\end{axis} 
\end{tikzpicture}
\caption{}
\label{fig:b_m_bm_easy}
\end{subfigure}
\begin{subfigure}[b]{0.45\textwidth}
\centering
\begin{tikzpicture}[scale=0.65]
\begin{axis}[ xlabel=Time (seconds), ylabel= \# problems solved,ymin=440, ymax=500, xmin=100, xmax=3600, legend pos=south east, legend cell align={left}, extra x ticks = {100}] 
\addplot[color=blue,line width=1] table[y=count, x=MIP_BM50, col sep=comma]{times.tex}; 
\addplot[color=green,line width=1,dashdotted] table[y=count, x=MIP_B50, col sep=comma]{times.tex}; 
\addplot[color=purple,line width=1,dashed] table[y=count, x=MIP_M50, col sep=comma]{times.tex};
\addplot[color=red,line width=1,dotted] table[y=count, x=BEN, col sep=comma]{times.tex}; 
\legend{\emph{FixB-FatM},\emph{FixB},\emph{FatM},\emph{Bench}}
\end{axis} 
\end{tikzpicture}
\caption{}
\label{fig:b_m_bm_hard}
\end{subfigure}
\caption{Computational burden of \emph{FixB}, \emph{FatM}, \emph{FixB-FatM} for $K$=50 and $\tau = 0$}
\label{fig:b_m_bm_results}
\end{figure}

It is also relevant to point out that the higher the value of $K$, the lower the chances of achieving unanimity on the status of switchable lines, and thus, the lower the number of binary variables that are fixed in the OTS problem. To illustrate this fact, Table \ref{tab:binaries} collects the results of approach \emph{FixB-FatM} for $\tau=0$ and for different values of $K$ including the average number of binary variables fixed to one or zero using the training data (\# bin). For $K=5$, 28 binary variables (out of the original 69 binary variables) are fixed in average, then leading to low computational times but a larger number of suboptimal instances. For $K=499$, only 8 binary variables are fixed (in average), no suboptimal solutions are obtained, but the computational time is increased. Figure \ref{fig:binaries} also illustrates the impact of $K$ on the performance of the \emph{FixB-FatM} approach. Note that setting $K$ equal to 5 yields the lowest computational times and all instances are solved in less than 100s. However, this method leads to 47 suboptimal solutions. On the other hand, if $K$ is set to 499, the maximum time reaches 400s but all instances are solved to optimality. 

\begin{table}[ht]
\setlength{\tabcolsep}{3pt}
\centering
\begin{tabular}{ccccccc}
\toprule 
 $K$ & \# opt & \# sub & gap-ave & gap-max & time (s) & \# bin \\
\midrule
 5  & 453 & 47 & 0.031 & 1.92 & 2.41 & 28.31 \\
 10  & 486 & 14 & 0.010 & 0.92 & 6.98 & 21.32 \\
 20  & 499 & 1 & 0.000 & 0.10 & 9.56 & 18.19 \\
 50  & 500 & 0 & - & - & 12.33 & 15.65 \\
 100 & 500 & 0 & - & - & 15.27 & 13.61 \\
 200 & 500 & 0 & - & - & 16.65 & 11.29 \\
 499 & 500 & 0 & - & - & 16.46 & 8.00 \\
 \bottomrule
\end{tabular}
\caption{Impact of $K$ on the performance of \emph{FixB-FatM} for $\tau = 0$}
\label{tab:binaries}
\end{table}

\begin{figure}[ht]
\captionsetup[subfigure]{labelformat=empty}
\centering
\begin{subfigure}[b]{0.45\textwidth}
\centering
\begin{tikzpicture}[scale=0.65]
\begin{axis}[ xlabel=Time (seconds), ylabel= \# problems solved,ymin=0, ymax=500, xmin=0, xmax=100, legend pos=south east, legend cell align={left}] 
\addplot[color=blue,line width=1] table[y=count, x=MIP_BM5, col sep=comma]{times.tex}; 
\addplot[color=green,line width=1,dashdotted] table[y=count, x=MIP_BM50, col sep=comma]{times.tex};
\addplot[color=purple,line width=1,dashed] table[y=count, x=MIP_BM500, col sep=comma]{times.tex};
\addplot[color=red,line width=1,dotted] table[y=count, x=BEN, col sep=comma]{times.tex};
\legend{$K=5$,$K=50$,$K=499$,\emph{Bench}}
\end{axis} 
\end{tikzpicture}
\caption{}
\label{fig:binaries_easy}
\end{subfigure}
\begin{subfigure}[b]{0.45\textwidth}
\centering
\begin{tikzpicture}[scale=0.65]
\begin{axis}[ xlabel=Time (seconds), ylabel= \# problems solved,ymin=440, ymax=500, xmin=100, xmax=500, legend pos=south east, legend cell align={left}] 
\addplot[color=blue,line width=1] table[y=count, x=MIP_BM5, col sep=comma]{times.tex}; 
\addplot[color=green,line width=1,dashdotted] table[y=count, x=MIP_BM50, col sep=comma]{times.tex};
\addplot[color=purple,line width=1,dashed] table[y=count, x=MIP_BM500, col sep=comma]{times.tex};
\addplot[color=red,line width=1,dotted] table[y=count, x=BEN, col sep=comma]{times.tex};
\legend{$K=5$,$K=50$,$K=499$,\emph{Bench}}
\end{axis} 
\end{tikzpicture}
\caption{}
\label{fig:binaries_hard}
\end{subfigure}
\caption{Impact of $K$ on the computational burden of \emph{FixB-FatM} for $\tau=0$}
\label{fig:binaries}
\end{figure}

While in the previous simulations $\tau$ was set to zero in all cases, increasing its value has the potential to fix a greater number of binary variables, thereby decreasing the time to solve the resulting OTS problem. This, however, comes at the cost of potentially increasing the number of infeasible and/or suboptimal instances. For $K=50$, Table \ref{tab:impact_tau} presents the simulation results for the \emph{FixB-FatM} method with various values of the threshold parameter $\tau$. The last column of this table ($\#$ bin) shows the average number of fixed binary variables, which logically rises with increasing values of $\tau$. However, the reduction in computational time is arguably minor and, certainly, may not justify the trade-off, as the gap values and the number of suboptimal instances increase significantly in contrast. Therefore, the relaxation of the unanimity condition in the proposed learning-based methods may not be deemed worthwhile.

\begin{table}[ht]
\setlength{\tabcolsep}{4pt}
\centering
\begin{tabular}{lccccccccc}
\toprule
&$K$ & $\tau$ & \# opt & \# sub & \# inf & gap-ave & gap-max & time (s) & \# bin \\
\midrule
\emph{FixB-FatM} &    50 &       0.00 &  500&   0&  0&     - &     - &   12.33 &  15.65 \\
\emph{FixB-FatM} &    50 &       0.01 &  500&   0&  0&     - &     - &   11.23 &  15.65 \\
\emph{FixB-FatM} &    50 &       0.02 &  497&   3&  0&     0.001 &   0.26 &   9.73 &  17.95 \\
\emph{FixB-FatM} &    50 &       0.05 &  493&   7&  0&     0.003 &    0.32 &   10.56 &  19.17 \\
\emph{FixB-FatM} &    50 &       0.10 &  486&  14&  0&     0.009 &     0.76 &   6.99 &  22.16 \\
\emph{FixB-FatM} &    50 &       0.20 &  454&  46&  0&     0.031 &     1.96 &   2.43 &  27.61 \\
\bottomrule
\end{tabular}
\caption{Impact of threshold $\tau$ on \emph{FixB-FatM} approach} 
\label{tab:impact_tau}
\end{table}

Next, we analyze the results of the two remaining approaches: the \emph{FixB-AngM} approach that uses the nearest neighbors to fix some binary variables and all the elements in the training to learn the big-M values as explained in Section~\ref{subs:proposed}, and the \emph{AngM} approach described in the same section. The results of these two methods for $\lambda=1$ are provided in Table \ref{tab:ang_results} and allow us to draw some interesting conclusions. First, both approaches lead to suboptimal solutions for some instances. This is understandable since, as explained in Section~\ref{subs:proposed}, these methods set the big-M constants fully relying on the maximum angle difference observed in the training set. Therefore, if the training set does not include an instance in which the actual maximum angle difference realizes, then the learned values of the big-Ms may leave the optimal solution out of the feasible region. In other words, while this strategy usually leads to very tight big-M values, it also increases the probability of having suboptimal or even infeasible solutions. This strategy is substantially different from approaches \emph{FatM} and \emph{FixB-FatM} that learn shorter paths of connected lines based on the optimal solution of the OTS problem for the training data and recompute the big-M constants using \eqref{eq: big M Fattahi}. Since shorter paths are only updated under the unanimity of the nearest neighbors, this strategy leads to more conservative big-M values and, consequently, larger feasibility regions and computational times. These facts are confirmed by comparing Tables \ref{tab:binaries} and \ref{tab:ang_results}. For instance, for $K=50$, \emph{FixB-FatM} solves all instances to optimality and takes 12.33s in average, the \emph{FixB-AngM} yields five suboptimal solutions but the average computational times is reduced to 0.7s only. The third relevant fact arises from the comparison of the average computational times of the two approaches in Table~\ref{tab:ang_results}. As observed, these times are particularly similar for all values of $K$. This leads us to conclude that the obtained big-M constants are so tight that fixing some binary variables does not have a significant impact on the computational burden. For completeness, Figure \ref{fig:ang_results} compares, for $\lambda=1$, the number of problems solved by \emph{FixB-AngM} for 50 neighbors and \emph{AngM} with the benchmark. Notice that  these two methodologies are able to solve most instances in less than 5 seconds, while only 250 instances are solved by the benchmark in that time. This figure also proves that fixing the binary variables has a negligible effect on the computational savings. 

\begin{table}[ht]
\setlength{\tabcolsep}{3pt}
\centering
\begin{tabular}{lcccccc}
\toprule
&$K$ & \# opt & \# sub & gap-ave & gap-max & time (s) \\
\midrule
\emph{FixB-AngM} & 5 & 450 & 50 & 0.033 & 1.92 & 0.41 \\
\emph{FixB-AngM} & 10 & 482 & 18 & 0.011& 0.92 & 0.59 \\
\emph{FixB-AngM} & 20 & 494 & 6 & 0.002 & 0.39 & 0.61 \\
\emph{FixB-AngM} & 50 & 495 & 5 & 0.002 & 0.39 & 0.70 \\
\emph{FixB-AngM} & 100 & 495 & 5 & 0.002 & 0.39 & 0.71\\
\emph{FixB-AngM} & 200 & 495 & 5 & 0.002 & 0.39 & 0.70\\
\emph{FixB-AngM} & 499 & 495 & 5 & 0.002 & 0.39 & 0.71\\
\emph{AngM} &-& 495 & 5 & 0.002 & 0.39 & 0.88\\
\bottomrule
\end{tabular}
\caption{Performance of approaches \emph{FixB-AngM} and \emph{AngM} for $\lambda=1$}
\label{tab:ang_results}
\end{table}

\begin{figure}[ht]
\centering
\begin{tikzpicture}[scale=0.65]
\begin{axis}[ xlabel=Time (seconds), ylabel= \# problems solved,ymin=0, ymax=500, xmin=0, xmax=5, legend pos=south east, legend cell align={left}] 
\addplot[color=blue,line width=1] table[y=count, x=ANG-BA50, col sep=comma]{times_small.tex};
\addplot[color=green,line width=1,dashed] table[y=count, x=ANG-K500, col sep=comma]{times_small.tex};
\addplot[color=red,line width=1,dotted] table[y=count, x=BEN, col sep=comma]{times_small.tex};
\legend{\emph{FixB-AngM},\emph{AngM},\emph{Bench}}
\end{axis} 
\end{tikzpicture}
\caption{Computational burden of \emph{FixB-AngM} and \emph{AngM} for $K=50$ and $\lambda=1$}
\label{fig:ang_results}
\end{figure}

To reduce the number of suboptimal instances, \emph{AngM}  can be run with values of the multiplying factor $\lambda$ higher than 1. Table \ref{tab:impact_lambda} compiles the simulation results for \emph{AngM} with various values of $\lambda$. It is observed that a slight increase in the big-M values above those learned from historical data has a minimal impact on computational time, but reduces the number of suboptimal instances. Remarkably, even for $\lambda=1.1$, all instances are solved optimally by \emph{AngM}.

\begin{table}[ht]
\setlength{\tabcolsep}{6pt}
\centering
\begin{tabular}{lccccccc}
\toprule
& $\lambda$ & \# opt & \# sub & \# inf & gap-ave & gap-max & time (s)  \\
\midrule
\emph{AngM} &       1.0 &  495&  5&  0&     0.002 &     0.39 &   0.88 \\
\emph{AngM} &       1.1 &  500&  0&  0&      - &      - &   0.78 \\
\emph{AngM} &       1.2 &  500&  0&  0&      - &      - &   0.82 \\
\emph{AngM} &       1.5 &  500&  0&  0&      - &      - &   1.21 \\
\bottomrule
\end{tabular}
\caption{Impact of factor $\lambda$ on \emph{AngM} approach} 
\label{tab:impact_lambda}
\end{table}

To further illustrate the performance of the two data-driven strategies to learn the big-M constants, Table \ref{tab:big-M} provides, for ten of the switchable lines, the big-M values for approaches \emph{Bench}, \emph{FixB-FatM} for $K=50,\tau=0$ and \emph{AngM} for $\lambda=1$. For the first two methods, $\underline{M}_{nm}$ and $\overline{M}_{nm}$ are symmetric for all lines, whereas approach \emph{AngM} computes asymmetric values as explained in Section \ref{subs:proposed}. Since the learned large constants may change for each instance, Table \ref{tab:big-M} includes value ranges. Thanks to the status of switchable lines of the nearest neighbors, the \emph{FixB-FatM} approach is able to reduce the shortest paths used in \eqref{eq: big M Fattahi} and significantly decrease the values of the big-Ms for some lines. For lines 2, 58 and 103, these values remain, however, unaltered. The approach \emph{AngM} learns from the observed angle differences and therefore, the big-M are tightened even further. In fact, for lines 58, 85, 135, 164, this methodology is able to infer the direction of the power flow through these lines and consequently one of the big-M values is set to 0. This bound reduction effectively tightens the DC-OTS model \eqref{eq:ots_mip} and significantly reduces its computational burden.

\begin{table}[ht]
\setlength{\tabcolsep}{6pt}
\centering
\begin{tabular}{ccccc}
\toprule
 line & \emph{Bench} & \emph{FixB-FatM} & \multicolumn{2}{c}{\emph{AngM}} \\
 & $-\underline{M}=\overline{M}$ &  $-\underline{M}=\overline{M}$ & $\overline{M}$ & $-\underline{M}$\\
\midrule
 2 & 1080 & 1080 & [212,218] & [388,383] \\
 23 & 10267 & [6615,10267] & [1441,1575] & [639,607] \\
 28 & 16806 & [7434,16806] & [553,628] & [604,510] \\
 31 & 1417 & [1309,1417] & [248,252] & [176,175] \\
 46 & 5279 & [2287,5279] & [289,325] & [34,9] \\
 58 & 247 & 247  & 0 & 81  \\
 85 & 776 & [0,776] & [376,391] & 0 \\
 103 & 486 & 486  & [184,185] & 381  \\
 135 & 1458 & [0,294] & [122,127] & 0 \\
 164 & 3231 & [0,837] & 0 & [115,114] \\
\bottomrule
\end{tabular}
\caption{Comparison of big-M values for \emph{Bench}, \emph{FixB-FatM}, \emph{AngM}}
\label{tab:big-M}
\end{table}

After this in-depth analysis of the simulation results for the \textit{Unif10} database, we can conclude that the most promising approaches are \textit{Linear} with $K=499$, \textit{FixB-FatM} with $K=50$ and $\tau=0$, and \textit{AngM} with $\lambda=1.1$. Table~\ref{tab:summary_unif10} summarizes the computational results of these approaches. The \textit{Linear} approach is the fastest, but returns 373 suboptimal instances, a maximum gap of 0.71\% and an average gap that is four times the target value of 0.01\%. On the other hand, \textit{FixB-FatM} and \textit{AngM} achieve the optimal solution for all instances. Besides, \textit{AngM} reports the lowest computational time, which is in fact slightly above that of the \textit{Linear} approach. 

\begin{table}[ht]
\setlength{\tabcolsep}{6pt}
\centering
\begin{tabular}{lccccccc}
\hline
&$K$ & \# opt & \# sub & \# inf  & gap-ave & gap-max & time (s) \\
\hline
\emph{Linear} &    499 &  127&   373&  0&      0.043&     0.71 &    0.39 \\
\emph{FixB-FatM}($\tau=0$) &    50 &  500&    0&  0&      - &  - &   12.33 \\
\emph{AngM}($\lambda=1.1$) &   - &  500&    0&  0&     - &  - & 0.78 \\
\hline
\end{tabular}
\caption{Summary of computational results for the \emph{Unif10} database} 
\label{tab:summary_unif10}
\end{table}

To conclude this section, we remark that the primary goal of these learning procedures is to swiftly generate solutions needed for \textit{online} applications. However, it is crucial to note that the rapid solutions obtained are not directly included in the training data. As new demand levels materialize over time, each instance must undergo an \textit{offline} optimization using the benchmark approach to achieve optimality before integrating its corresponding solution into the expanding training set.

\subsection{Impact of demand variability and correlation}\label{sub:impact_variability}

As mentioned earlier, the variability and correlation of nodal demand levels can influence the performance of the learning-based methods compared in this paper. Specifically, an increase in demand variability relative to nominal values is expected to reduce the accuracy of any learning method, given the same size of the training dataset. Conversely, a higher correlation among demand levels at different nodes in the network simplifies the learning task, thanks to a more pronounced data structure.

Table \ref{tab:unif20} compiles the simulation results of various methods for the \emph{Unif20} database, which has a higher variability than  \emph{Unif10}. While none of the methods applied to the \emph{Unif10} database result in any infeasible instances, this is not the case for the \emph{Unif20} dataset. The fifth column of the table indicates the number of infeasible instances for each approach. It is worth noting that, for this dataset, the benchmark fails to achieve optimality for 43 instances within one hour, resulting in an average mip-gap of 0.50\% and a maximum mip-gap of 2.40\%. The average time required by the benchmark method is 510.9s. For the \emph{Unif20} database, which includes the most challenging instances, we do observe a few cases where some of the learning-based methods produce slightly improved integer solutions compared to \emph{Bench}. However, for consistency, the reported gaps in this case study are calculated using the solutions identified by the \textit{Bench} approach as optimal. The simulation results in Table \ref{tab:unif20} yield noteworthy observations. Firstly, as anticipated, increasing the variability of demand levels leads to a rise in the number of suboptimal and infeasible instances. For instance,  \emph{FixB-FatM} with $K=5$ produced 47 suboptimal instances for the \emph{Unif10} database. However, for the \emph{Unif20} database, this method resulted in 153 suboptimal instances and 4 infeasible problems. The maximum gap for this approach has also increased from 1.92\% to 5.93\%. Secondly, augmenting the number of closest neighbors diminishes the number of infeasible instances, as binary variables are fixed only under the unanimity condition. Indeed, the \emph{FixB-AngM} approach exhibits no infeasible instances when $K$ is increased from 5 to 50. Similarly, the Linear approach avoids any infeasible instance for a value of $K=499$. This suggests that these approaches are not particularly suitable for high variability in parameters or a low number of training instances. Thirdly, the \emph{Linear} approach is very fast, but involves large average and maximum gap values, even for $K=499$. Finally, considering both the number of suboptimal instances, the average and maximum gaps, and the average computational time, it can be concluded that the \emph{AngM} with $\lambda=1.1$ method exhibits superior performance for the \emph{Unif20} database.

\begin{table}[ht]
\setlength{\tabcolsep}{6pt}
\centering
\begin{tabular}{lccccccc}
\toprule
&$K$ & \# opt & \# sub & \# inf & gap-ave & gap-max & time (s) \\
\midrule
\emph{Linear} &     5 &    2 &  490&   8&     0.923 &     6.27 &    0.00 \\
\emph{Linear} &    50 &   18&  481&   1&     0.276 &     2.63 &    0.04 \\
\emph{Linear} &    499 &   44&  456&   0&     0.153 &     1.50 &    0.33 \\
\emph{FixB-FatM}($\tau=0$) &     5 &  343 &  153&   4&     0.179 &     5.93 &    4.99 \\
\emph{FixB-FatM}($\tau=0$) &    50 &  496&    4&   0&     0.002 &     0.72 &  115.02 \\
\emph{FixB-FatM}($\tau=0$) &    499 &  499&    1&   0&     0.000 &     0.02 &  105.54 \\
\emph{AngM}($\lambda=1$) &   - &  492&    8&    0&     0.008 &     2.00 &    1.38 \\
\emph{AngM}($\lambda=1.1$) &   - &  499&    1&    0&     0.000 &     0.02 &    1.93 \\
\bottomrule
\end{tabular}
\caption{Computational results for the \emph{Unif20} database} 
\label{tab:unif20}
\end{table}

Despite the insightful findings presented in Table \ref{tab:unif20}, one could argue that electricity demand in real power systems exhibits a higher spatial correlation. Therefore, utilizing uncorrelated probability distributions for the nodal demands may diverge from reality. To address this concern, in Table  \ref{tab:corr} we present results analogous to those in Table \ref{tab:unif20} where demand levels are randomly sampled from a multinormal distribution with a correlation matrix computed using data from \cite{joswig2021opflearndata}. For the \emph{Normal} dataset, the benchmark approach fails to solve 30 instances within one hour, yielding an average mip-gap of 0.49\% and a maximum mip-gap of 2.45\%. Besides, the benchmark approach takes an average time of 289.5 seconds. As with the \textit{Unif10} database, none of the learning-based approaches improves the solution found by the \emph{Bench} approach in one hour for any of the 500 instances of the \textit{Normal} database. In this more realistic setting, we observe that there are no infeasible instances for any of the methods, while most methods result in some suboptimal instances. Notably, the computational times required by \emph{Linear}, \emph{FixB-FatM} and \emph{AngM} are of the same order of magnitude. However, while the \emph{Linear} approach returns  suboptimal instances for the three values of $K$, the proposed methodologies \emph{FixB-FatM} with $K=499$ and $\tau=0$, and \emph{AngM} with $\lambda=1.1$ are able to solve the 500 instances to global optimality. This underscores the efficacy of learning-based procedures in delivering rapid solutions that closely approximate the original solution for the OTS problem, even with realistic demand level variability.

\begin{table}[ht]
\setlength{\tabcolsep}{6pt}
\centering
\begin{tabular}{lccccccc}
\toprule
&$K$ & \# opt & \# sub & \# inf & gap-ave & gap-max & time (s) \\
\midrule
\emph{Linear} &     5 &  164&  336&  0&     0.024 &     0.47 &    0.00 \\
\emph{Linear} &    50 &  446&   54&  0&     0.004 &     0.37 &    0.04 \\
\emph{Linear} &    499 &  488&   12&  0&     0.001 &     0.11 &    0.41 \\
\emph{FixB-FatM}($\tau=0$) &     5 &  493&    7&  0&     0.003 &     0.46 &    0.30 \\
\emph{FixB-FatM}($\tau=0$) &    50 &  499&    1&  0&     0.000 &     0.17 &    0.57 \\
\emph{FixB-FatM}($\tau=0$) &    499 &  500&    0&  0&     - &     - &    0.56 \\
\emph{AngM}($\lambda=1$) &   - &  495&    5&  0&     0.002 &     0.51 &    0.26 \\
\emph{AngM}($\lambda=1.1$) &   - &  500&    0&  0&     - &     - &    0.29 \\
\bottomrule
\end{tabular}
\caption{Computational results for the \emph{Normal} database} 
\label{tab:corr}
\end{table}

\section{Conclusions and further research} \label{sec:concl}
In the field of power systems, the optimal transmission switching problem (OTS) determines the on/off status of transmission lines to reduce the operating cost. The OTS problem can be formulated as a mixed-integer linear program (MILP)  that includes large enough constants. This problem belongs to the NP-hard class and its computational burden is, consequently, significant even for small networks. While \textit{pure} end-to-end learning approaches can solve the OTS problem extremely fast, the obtained solutions are usually suboptimal, or even infeasible. Alternatively, we propose in this paper some learning-based approaches that reduce the computational burden of the MILP model by leveraging information of previously solved instances. These computational savings arise from the fact that some binary variables are fixed and tighter big-M values are found. Numerical simulations on a 118-bus power network show that the first proposed approach is able to solve all instances to optimality in less than 300 seconds, while the benchmark approach is unable to solve all of them in 3600 seconds. The second approach we propose is more \textit{aggressive} and solves all instances in less than 10 seconds, but 1\% of them do not reach the optimal solution. We also assess the performance of the proposed learning-based approaches under increased demand variability and correlation.

All the learning approaches presented in this paper utilize the $K$nn algorithm and the $l_2$ norm distance. The exploration of different machine learning methods and/or distances is left as a potential avenue for future research. In this paper, we introduce a machine learning approach that leverages the structural patterns observed in past DC-OTS instances to improve the performance of new problems. However, the solver hyperparameters are set to default values. Future research could explore utilizing the data information not only to exploit the problem structure but also to finely tune solver hyperparameters, as demonstrated in \cite{lodi2017learning, cappart2023combinatorial}. Additionally, our study assumes the use of DC approximations for power flow equations. A potential research direction involves addressing the more challenging AC-OTS problem, considering data-driven strategies to simplify it into a DC-OTS format, akin to approaches presented in \cite{parmentier2022learning}.

\bmhead{Acknowledgments}

This work was supported in part by the European Research Council (ERC) under the EU Horizon 2020 research and innovation program (grant agreement No. 755705), in part by the Spanish Ministry of Science and Innovation (AEI/10.13039/501100011033) through project PID2020-115460GB-I00, and in part by the Research Program for Young Talented Researchers of the University of Málaga under Project B1-2020-15. Finally, the authors thankfully acknowledge the computer resources, technical expertise, and assistance provided by the SCBI (Supercomputing and Bioinformatics) center of the University of Málaga.

\bmhead{Conflict of interest} The authors declare no conflict of interest.

\bibliography{sn-bibliography}


\begin{thebibliography}{29}
\ifx \bisbn   \undefined \def \bisbn  #1{ISBN #1}\fi
\ifx \binits  \undefined \def \binits#1{#1}\fi
\ifx \bauthor  \undefined \def \bauthor#1{#1}\fi
\ifx \batitle  \undefined \def \batitle#1{#1}\fi
\ifx \bjtitle  \undefined \def \bjtitle#1{#1}\fi
\ifx \bvolume  \undefined \def \bvolume#1{\textbf{#1}}\fi
\ifx \byear  \undefined \def \byear#1{#1}\fi
\ifx \bissue  \undefined \def \bissue#1{#1}\fi
\ifx \bfpage  \undefined \def \bfpage#1{#1}\fi
\ifx \blpage  \undefined \def \blpage #1{#1}\fi
\ifx \burl  \undefined \def \burl#1{\textsf{#1}}\fi
\ifx \doiurl  \undefined \def \doiurl#1{\url{https://doi.org/#1}}\fi
\ifx \betal  \undefined \def \betal{\textit{et al.}}\fi
\ifx \binstitute  \undefined \def \binstitute#1{#1}\fi
\ifx \binstitutionaled  \undefined \def \binstitutionaled#1{#1}\fi
\ifx \bctitle  \undefined \def \bctitle#1{#1}\fi
\ifx \beditor  \undefined \def \beditor#1{#1}\fi
\ifx \bpublisher  \undefined \def \bpublisher#1{#1}\fi
\ifx \bbtitle  \undefined \def \bbtitle#1{#1}\fi
\ifx \bedition  \undefined \def \bedition#1{#1}\fi
\ifx \bseriesno  \undefined \def \bseriesno#1{#1}\fi
\ifx \blocation  \undefined \def \blocation#1{#1}\fi
\ifx \bsertitle  \undefined \def \bsertitle#1{#1}\fi
\ifx \bsnm \undefined \def \bsnm#1{#1}\fi
\ifx \bsuffix \undefined \def \bsuffix#1{#1}\fi
\ifx \bparticle \undefined \def \bparticle#1{#1}\fi
\ifx \barticle \undefined \def \barticle#1{#1}\fi
\bibcommenthead
\ifx \bconfdate \undefined \def \bconfdate #1{#1}\fi
\ifx \botherref \undefined \def \botherref #1{#1}\fi
\ifx \url \undefined \def \url#1{\textsf{#1}}\fi
\ifx \bchapter \undefined \def \bchapter#1{#1}\fi
\ifx \bbook \undefined \def \bbook#1{#1}\fi
\ifx \bcomment \undefined \def \bcomment#1{#1}\fi
\ifx \oauthor \undefined \def \oauthor#1{#1}\fi
\ifx \citeauthoryear \undefined \def \citeauthoryear#1{#1}\fi
\ifx \endbibitem  \undefined \def \endbibitem {}\fi
\ifx \bconflocation  \undefined \def \bconflocation#1{#1}\fi
\ifx \arxivurl  \undefined \def \arxivurl#1{\textsf{#1}}\fi
\csname PreBibitemsHook\endcsname

\bibitem{o2005dispatchable}
\begin{barticle}
\bauthor{\bsnm{O'Neill}, \binits{R.P.}},
\bauthor{\bsnm{Baldick}, \binits{R.}},
\bauthor{\bsnm{Helman}, \binits{U.}},
\bauthor{\bsnm{Rothkopf}, \binits{M.H.}},
\bauthor{\bsnm{Stewart}, \binits{W.}}:
\batitle{Dispatchable transmission in {RTO} markets}.
\bjtitle{IEEE Transactions on Power Systems}
\bvolume{20}(\bissue{1}),
\bfpage{171}--\blpage{179}
(\byear{2005})
\end{barticle}
\endbibitem

\bibitem{fisher2008optimal}
\begin{barticle}
\bauthor{\bsnm{Fisher}, \binits{E.B.}},
\bauthor{\bsnm{O'Neill}, \binits{R.P.}},
\bauthor{\bsnm{Ferris}, \binits{M.C.}}:
\batitle{Optimal transmission switching}.
\bjtitle{IEEE Transactions on Power Systems}
\bvolume{23}(\bissue{3}),
\bfpage{1346}--\blpage{1355}
(\byear{2008})
\end{barticle}
\endbibitem

\bibitem{kocuk2016cycle}
\begin{barticle}
\bauthor{\bsnm{Kocuk}, \binits{B.}},
\bauthor{\bsnm{Jeon}, \binits{H.}},
\bauthor{\bsnm{Dey}, \binits{S.S.}},
\bauthor{\bsnm{Linderoth}, \binits{J.}},
\bauthor{\bsnm{Luedtke}, \binits{J.}},
\bauthor{\bsnm{Sun}, \binits{X.A.}}:
\batitle{A cycle-based formulation and valid inequalities for {DC} power
  transmission problems with switching}.
\bjtitle{Operations Research}
\bvolume{64}(\bissue{4}),
\bfpage{922}--\blpage{938}
(\byear{2016})
\end{barticle}
\endbibitem

\bibitem{fattahi2019bound}
\begin{barticle}
\bauthor{\bsnm{Fattahi}, \binits{S.}},
\bauthor{\bsnm{Lavaei}, \binits{J.}},
\bauthor{\bsnm{Atamtürk}, \binits{A.}}:
\batitle{A bound strengthening method for optimal transmission switching in
  power systems}.
\bjtitle{IEEE Transactions on Power Systems}
\bvolume{34}(\bissue{1}),
\bfpage{280}--\blpage{291}
(\byear{2019})
\end{barticle}
\endbibitem

\bibitem{ruiz2016security}
\begin{barticle}
\bauthor{\bsnm{Ruiz}, \binits{P.A.}},
\bauthor{\bsnm{Goldis}, \binits{E.}},
\bauthor{\bsnm{Rudkevich}, \binits{A.M.}},
\bauthor{\bsnm{Caramanis}, \binits{M.C.}},
\bauthor{\bsnm{Philbrick}, \binits{C.R.}},
\bauthor{\bsnm{Foster}, \binits{J.M.}}:
\batitle{Security-constrained transmission topology control milp formulation
  using sensitivity factors}.
\bjtitle{IEEE Transactions on Power Systems}
\bvolume{32}(\bissue{2}),
\bfpage{1597}--\blpage{1605}
(\byear{2016})
\end{barticle}
\endbibitem

\bibitem{dey2022node}
\begin{barticle}
\bauthor{\bsnm{Dey}, \binits{S.S.}},
\bauthor{\bsnm{Kocuk}, \binits{B.}},
\bauthor{\bsnm{Redder}, \binits{N.}}:
\batitle{Node-based valid inequalities for the optimal transmission switching
  problem}.
\bjtitle{Discrete Optimization}
\bvolume{43},
\bfpage{100683}
(\byear{2022})
\end{barticle}
\endbibitem

\bibitem{liu2012heuristic}
\begin{barticle}
\bauthor{\bsnm{Liu}, \binits{C.}},
\bauthor{\bsnm{Wang}, \binits{J.}},
\bauthor{\bsnm{Ostrowski}, \binits{J.}}:
\batitle{Heuristic prescreening switchable branches in optimal transmission
  switching}.
\bjtitle{IEEE Transactions on Power Systems}
\bvolume{27}(\bissue{4}),
\bfpage{2289}--\blpage{2290}
(\byear{2012})
\end{barticle}
\endbibitem

\bibitem{barrows2012computationally}
\begin{bchapter}
\bauthor{\bsnm{Barrows}, \binits{C.}},
\bauthor{\bsnm{Blumsack}, \binits{S.}},
\bauthor{\bsnm{Bent}, \binits{R.}}:
\bctitle{Computationally efficient optimal transmission switching: Solution
  space reduction}.
In: \bbtitle{2012 IEEE Power and Energy Society General Meeting},
pp. \bfpage{1}--\blpage{8}
(\byear{2012}).
\bcomment{IEEE}
\end{bchapter}
\endbibitem

\bibitem{flores2020alternative}
\begin{barticle}
\bauthor{\bsnm{Flores}, \binits{M.}},
\bauthor{\bsnm{Macedo}, \binits{L.H.}},
\bauthor{\bsnm{Romero}, \binits{R.}}:
\batitle{Alternative mathematical models for the optimal transmission switching
  problem}.
\bjtitle{IEEE Systems Journal}
\bvolume{15}(\bissue{1}),
\bfpage{1245}--\blpage{1255}
(\byear{2020})
\end{barticle}
\endbibitem

\bibitem{fuller2012fast}
\begin{barticle}
\bauthor{\bsnm{Fuller}, \binits{J.D.}},
\bauthor{\bsnm{Ramasra}, \binits{R.}},
\bauthor{\bsnm{Cha}, \binits{A.}}:
\batitle{Fast heuristics for transmission-line switching}.
\bjtitle{IEEE Transactions on Power Systems}
\bvolume{27}(\bissue{3}),
\bfpage{1377}--\blpage{1386}
(\byear{2012})
\end{barticle}
\endbibitem

\bibitem{crozier2022feasible}
\begin{barticle}
\bauthor{\bsnm{Crozier}, \binits{C.}},
\bauthor{\bsnm{Baker}, \binits{K.}},
\bauthor{\bsnm{Toomey}, \binits{B.}}:
\batitle{Feasible region-based heuristics for optimal transmission switching}.
\bjtitle{Sustainable Energy, Grids and Networks}
\bvolume{30},
\bfpage{100628}
(\byear{2022})
\end{barticle}
\endbibitem

\bibitem{hinneck2022optimal}
\begin{botherref}
\oauthor{\bsnm{Hinneck}, \binits{A.}},
\oauthor{\bsnm{Pozo}, \binits{D.}}:
Optimal transmission switching: improving exact algorithms by parallel
  incumbent solution generation.
IEEE Transactions on Power Systems
(2022)
\end{botherref}
\endbibitem

\bibitem{johnson2020knearest}
\begin{botherref}
\oauthor{\bsnm{Johnson}, \binits{E.S.}},
\oauthor{\bsnm{Ahmed}, \binits{S.}},
\oauthor{\bsnm{Dey}, \binits{S.S.}},
\oauthor{\bsnm{Watson}, \binits{J.-P.}}:
A K-Nearest Neighbor Heuristic for Real-Time DC Optimal Transmission Switching.
arXiv
(2021).
\doiurl{10.48550/ARXIV.2003.10565}.
\url{https://arxiv.org/abs/2003.10565}
\end{botherref}
\endbibitem

\bibitem{yang2019line}
\begin{bchapter}
\bauthor{\bsnm{Yang}, \binits{Z.}},
\bauthor{\bsnm{Oren}, \binits{S.}}:
\bctitle{Line selection and algorithm selection for transmission switching by
  machine learning methods}.
In: \bbtitle{2019 IEEE Milan PowerTech},
pp. \bfpage{1}--\blpage{6}
(\byear{2019}).
\bcomment{IEEE}
\end{bchapter}
\endbibitem

\bibitem{han2022learning}
\begin{botherref}
\oauthor{\bsnm{Han}, \binits{T.}},
\oauthor{\bsnm{Hill}, \binits{D.}}:
Learning-based topology optimization of power networks.
IEEE Transactions on Power Systems
(2022)
\end{botherref}
\endbibitem

\bibitem{bugaje2023real}
\begin{barticle}
\bauthor{\bsnm{Bugaje}, \binits{A.-A.B.}},
\bauthor{\bsnm{Cremer}, \binits{J.L.}},
\bauthor{\bsnm{Strbac}, \binits{G.}}:
\batitle{Real-time transmission switching with neural networks}.
\bjtitle{IET Generation, Transmission \& Distribution}
\bvolume{17}(\bissue{3}),
\bfpage{696}--\blpage{705}
(\byear{2023})
{\href{https://arxiv.org/abs/https://ietresearch.onlinelibrary.wiley.com/doi/pdf/10.1049/gtd2.12698}{{https://ietresearch.onlinelibrary.wiley.com/doi/pdf/10.1049/gtd2.12698}}}.
\doiurl{10.1049/gtd2.12698}
\end{barticle}
\endbibitem

\bibitem{bengio2021machine}
\begin{barticle}
\bauthor{\bsnm{Bengio}, \binits{Y.}},
\bauthor{\bsnm{Lodi}, \binits{A.}},
\bauthor{\bsnm{Prouvost}, \binits{A.}}:
\batitle{Machine learning for combinatorial optimization: A methodological tour
  d’horizon}.
\bjtitle{European Journal of Operational Research}
\bvolume{290}(\bissue{2}),
\bfpage{405}--\blpage{421}
(\byear{2021}).
\doiurl{10.1016/j.ejor.2020.07.063}
\end{barticle}
\endbibitem

\bibitem{parmentier2022learning}
\begin{barticle}
\bauthor{\bsnm{Parmentier}, \binits{A.}}:
\batitle{Learning to approximate industrial problems by operations research
  classic problems}.
\bjtitle{Operations Research}
\bvolume{70}(\bissue{1}),
\bfpage{606}--\blpage{623}
(\byear{2022})
\end{barticle}
\endbibitem

\bibitem{pineda2020data}
\begin{barticle}
\bauthor{\bsnm{Pineda}, \binits{S.}},
\bauthor{\bsnm{Morales}, \binits{J.M.}},
\bauthor{\bsnm{Jiménez-Cordero}, \binits{A.}}:
\batitle{Data-driven screening of network constraints for unit commitment}.
\bjtitle{IEEE Transactions on Power Systems}
\bvolume{35}(\bissue{5}),
\bfpage{3695}--\blpage{3705}
(\byear{2020}).
\doiurl{10.1109/TPWRS.2020.2980212}
\end{barticle}
\endbibitem

\bibitem{jimenez2022warm}
\begin{barticle}
\bauthor{\bsnm{Jiménez-Cordero}, \binits{A.}},
\bauthor{\bsnm{Morales}, \binits{J.M.}},
\bauthor{\bsnm{Pineda}, \binits{S.}}:
\batitle{Warm-starting constraint generation for mixed-integer optimization: A
  machine learning approach}.
\bjtitle{Knowledge-Based Systems}
\bvolume{253},
\bfpage{109570}
(\byear{2022}).
\doiurl{10.1016/j.knosys.2022.109570}
\end{barticle}
\endbibitem

\bibitem{hedman2012flexible}
\begin{botherref}
\oauthor{\bsnm{Hedman}, \binits{K.W.}},
\oauthor{\bsnm{Oren}, \binits{S.S.}},
\oauthor{\bsnm{O’Neill}, \binits{R.P.}}:
Flexible transmission in the smart grid: optimal transmission switching.
Handbook of networks in power systems I,
523--553
(2012)
\end{botherref}
\endbibitem

\bibitem{cormen2022introduction}
\begin{bbook}
\bauthor{\bsnm{Cormen}, \binits{T.H.}},
\bauthor{\bsnm{Leiserson}, \binits{C.E.}},
\bauthor{\bsnm{Rivest}, \binits{R.L.}},
\bauthor{\bsnm{Stein}, \binits{C.}}:
\bbtitle{Introduction to Algorithms}.
\bpublisher{MIT press},
\blocation{Cambridge, Massachusetts}
(\byear{2022})
\end{bbook}
\endbibitem

\bibitem{hastie2009elements}
\begin{bbook}
\bauthor{\bsnm{Hastie}, \binits{T.}},
\bauthor{\bsnm{Tibshirani}, \binits{R.}},
\bauthor{\bsnm{Friedman}, \binits{J.H.}},
\bauthor{\bsnm{Friedman}, \binits{J.H.}}:
\bbtitle{The Elements of Statistical Learning: Data Mining, Inference, and
  Prediction}
vol. \bseriesno{2}.
\bpublisher{Springer},
\blocation{New Yourk}
(\byear{2009})
\end{bbook}
\endbibitem

\bibitem{blumsack2006network}
\begin{bbook}
\bauthor{\bsnm{Blumsack}, \binits{S.}}:
\bbtitle{Network Topologies and Transmission Investment Under Electric-industry
  Restructuring}.
\bpublisher{Carnegie Mellon University},
\blocation{Pittsburgh, Pennsylvania}
(\byear{2006})
\end{bbook}
\endbibitem

\bibitem{joswig2021opflearndata}
\begin{botherref}
\oauthor{\bsnm{Joswig-Jones}, \binits{T.}},
\oauthor{\bsnm{Zamzam}, \binits{A.}},
\oauthor{\bsnm{Baker}, \binits{K.}}:
{OPFL}earndata: Dataset for learning {AC} optimal power flow.
Technical report,
NREL Data Catalog. Golden, CO: National Renewable Energy Laboratory
(2021)
\end{botherref}
\endbibitem

\bibitem{OASYS2023learning}
\begin{botherref}
\oauthor{\bsnm{OASYS}}:
{Learning\_Assisted\_Optimization\_for\_Transmission\_Switching}.
\url{https://github.com/groupoasys/Learning_Assisted_Optimization_for_Transmission_Switching}
(2023)
\end{botherref}
\endbibitem

\bibitem{gurobi}
\begin{botherref}
\oauthor{\bsnm{{Gurobi Optimization, LLC}}}:
{Gurobi Optimizer Reference Manual}
(2022).
\url{https://www.gurobi.com}
\end{botherref}
\endbibitem

\bibitem{lodi2017learning}
\begin{barticle}
\bauthor{\bsnm{Lodi}, \binits{A.}},
\bauthor{\bsnm{Zarpellon}, \binits{G.}}:
\batitle{On learning and branching: a survey}.
\bjtitle{Top}
\bvolume{25},
\bfpage{207}--\blpage{236}
(\byear{2017})
\end{barticle}
\endbibitem

\bibitem{cappart2023combinatorial}
\begin{barticle}
\bauthor{\bsnm{Cappart}, \binits{Q.}},
\bauthor{\bsnm{Ch{\'e}telat}, \binits{D.}},
\bauthor{\bsnm{Khalil}, \binits{E.B.}},
\bauthor{\bsnm{Lodi}, \binits{A.}},
\bauthor{\bsnm{Morris}, \binits{C.}},
\bauthor{\bsnm{Veli{\v{c}}kovi{\'c}}, \binits{P.}}:
\batitle{Combinatorial optimization and reasoning with graph neural networks}.
\bjtitle{Journal of Machine Learning Research}
\bvolume{24}(\bissue{130}),
\bfpage{1}--\blpage{61}
(\byear{2023})
\end{barticle}
\endbibitem

\end{thebibliography}

\end{document}